\documentclass{article}
  \usepackage{amsfonts,amssymb,mathrsfs,pb-diagram,makeidx}
  \usepackage[all]{xy}

  \makeindex

\def \1{{\bf 1}}

\def \al{\alpha}
\def \A{{\mathbb A}}
\def \Ad{{\rm Ad}}

\def \be{\beta}
\def \bs{\backslash}

\def \card{{\rm card}}

\def \C{{\mathbb C}}

\def \CCC{{\cal C}}
\def \CCD{{\cal D}}
\def \CD{{\mathscr D}}
\def \CE{{\mathscr E}}

\def \CF{{\mathscr F}}
\def \CG{{\mathscr G}}

\def \CCO{{\mathscr O}}
\def \CO{{\cal O}}
\def \cos{{\rm cos \hspace{2pt}}}

\def \d{\delta}
\def \D{\Delta}
\def \det{{\rm det \hspace{2pt}}}
\def \diag{{\rm diag}}

\def \ell{{\rm ell}}
\def \ep{\varepsilon}

\def \fin{{\rm fin}}
\def \F{{\mathbb F}}
\def \CF{{\mathfrak F}}

\def \ga{\gamma}
\def \Ga{\Gamma}

\def \GL{{\rm GL}}

\def \H{{\mathbb H}}

\def \hra{\hookrightarrow}

\def \ind{{\rm ind \hspace{2pt}}}
\def \Ind{{\rm Ind \hspace{2pt}}}
\def \Inv{{\rm Inv \hspace{2pt}}}

\def \ka{\kappa}

\def \la{\lambda}

\def \li{{\rm li}}

\def \Lin{{\rm Lin}}

\def \Mat{{\rm Mat}}
\def \mod{{\rm mod}}
\def \n{{{\mathfrak n}}}

\def \N{\mathbb N}

\def \ox{\otimes}
\def \p{{{\mathfrak p}}}

\def \prf{\noindent{\bf Proof: }}

\def \PGL{{\rm PGL}}

\def \qed{\ifmmode\eqno $\square$ \else\noproof\vskip 12pt plus 3pt minus 9pt \fi}
 \def\noproof{{\unskip\nobreak\hfill\penalty50\hskip2em\hbox{}%
     \nobreak\hfill $\square$ \parfillskip=0pt%
     \finalhyphendemerits=0\par}}
\def \Q{\mathbb Q}

\def \R{{\mathbb R}}

\def \ra{\rightarrow}

\def \reg{{\rm reg}}

\def \rw{{\rm rw}}
\def \rnw{{\rm rnw}}
\def \si{\sigma}
\def \Si{\Sigma}
\def \sin{{\rm sin \hspace{2pt}}}
\def \SL{{\rm SL}}

\def \SO{{\rm SO}}

\def \th{\theta}
\def \tr{{\hspace{1pt}\rm tr\hspace{2pt}}}
\def \Upp{{\rm Upp}}
\def \ut{\tilde{u}}

\def \x{\times}

\def \ze{\zeta}
\def \Z{\mathbb Z}
\def \={\ =\ }

\newcommand{\rez}[1]{\frac{1}{#1}}

\renewcommand{\matrix}[4]{\left( \begin{array}{cc}#1 & #2 \\ #3 & #4 \end{array}
            \right)}
\newcommand{\matrixfour}[4]
	   {\left( \begin{array}{cccc} 
	       #1 & \  & \  & \  \\ 
	       \  & #2 & \  & \  \\  
	       \  & \  & #3 & \  \\  
	       \  & \  & \  & #4 \end{array}\right)}
\newcommand{\matrixtwo}[2]{\matrix {#1}{}{}{#2}}
\newtheorem{theorem}{Theorem}[section]
\newtheorem{conjecture}[theorem]{Conjecture}
\newtheorem{lemma}[theorem]{Lemma}

\newtheorem{proposition}[theorem]{Proposition}

\begin{document}

 \title{Class Numbers of Orders in Quartic Fields}
 \author{Anton Deitmar \& Mark Pavey}
 \date{}
 \maketitle
 
\section*{Introduction}
Let $\CD$ be the set of all natural numbers $D\equiv 0,1\ \mod\,4$ with $D$ not a square.  Then $\CD$ is the set of all discriminants of orders in real quadratic fields.  For $D\in\CD$ the set
$$
\CO_D=\left\{\frac{x+y\sqrt{D}}{2}:x\equiv yD\ \mod\,2\right\}
$$
is an order in the real quadratic field $\Q(\sqrt{D})$ with discriminant $D$.  As $D$ varies, $\CO_D$ runs through the set of all orders of real quadratic fields.  For $D\in\CD$ let $h(\CO_D)$ denote the class number and $R(\CO_D)$ the regulator of the order $\CO_D$.  It was conjectured by Gauss (\cite{Gauss86}) and proved by Siegel (\cite{Siegel44}) that, as $x$ tends to infinity
$$
\sum_{{D\in\CD}\atop{D\leq x}}h(\CO_D)R(\CO_D)=\frac{\pi^2 x^{3/2}}{18\ze(3)}+O(x\log x),
$$
where $\ze$ is the Riemann zeta function.

For a long time it was believed to be impossible to separate the class number and the regulator in the summation.  However, in \cite{Sarnak82}, Theorem 3.1, Sarnak showed, using the Selberg trace formula, that as $x\ra\infty$ we have
$$
\sum_{{D\in\CD}\atop{R(\CO_D)\leq x}}h(\CO_D)\sim\frac{e^{2x}}{2x}.
$$
More sharply,
$$
\sum_{{D\in\CD}\atop{R(\CO_D)\leq x}}h(\CO_D)=L(2x)+O\left(e^{3x/2}x^2\right)
$$
as $x\ra\infty$, where $L(x)$ is the function
$
L(x)=\int_1^x\frac{e^t}{t}\,dt.
$
Sarnak established this result by identifying the regulators with lengths of closed geodesics of the modular curve $\SL_2(\Z)\bs\H$, where $\H$ denotes the upper half-plane, and using the Prime Geodesic Theorem for this Riemannian surface.  Actually Sarnak proved not this result but the analogue where $h(\CO)$ is replaced by the class number in the narrower sense and $R(\CO)$ by a ``regulator in the narrower sense".  But in Sarnak's proof the group $\SL_2(\Z)$ can be replaced by $\PGL_2(\Z)$ giving the above result.

The Prime Geodesic Theorem in this context gives an asymptotic formula for the number of closed geodesics on the surface $\SL_2(\Z)\bs\H$ with length less than or equal to $x>0$.  This formula is analogous to the asymptotic formula for the number of primes less than $x$ given in the Prime Number Theorem.  The Selberg zeta function (see \cite{Selberg56}) is used in the proof of the Prime Geodesic Theorem in a way analogous to the way the Riemann zeta function is used in the proof of the Prime Number Thoerem (see \cite{Chandrasekharan68}).  The required properties of the Selberg zeta function are deduced from the Selberg trace formula (\cite{Selberg56}).

It seems that following Sarnak's result no asymptotic results for class numbers in fields of degree greater than two were proven until in \cite{Deitmar02}, Theorem 1.1, the first named author proved an asymptotic formula for class numbers of orders in complex cubic fields, that is, cubic fields with one real embedding and one pair of complex conjugate embeddings.  This result can be stated as follows.

Let $S$ be a finite set of prime numbers containing at least two elements and let $C(S)$ be the set of all complex cubic fields $F$ such that all primes $p\in S$ are non-decomposed in $F$.  For $F\in C(S)$ let $O_F(S)$ be the set of all isomorphism classes of orders in $F$ which are maximal at all $p\in S$, ie. are such that the completion $\CO_p=\CO\ox\Z_p$ is the maximal order of the field $F_p=F\ox\Q_p$ for all $p\in S$.  Let $O(S)$ be the union of all $O_F(S)$, where $F$ ranges over $C(S)$.  For a field $F\in C(S)$ and an order $\CO\in O_F(S)$ define
$$
\la_S(\CO)=\la_S(F)=\prod_{p\in S}f_p(F),
$$
where $f_p(F)$ is the inertia degree of $p$ in $F$.  Let $R(\CO)$ denote the regulator and $h(\CO)$ the class number of the order $\CO$.

For $x>0$ we define
$$
\pi_S(x)=\sum_{{\CO\in O(S)}\atop{R(\CO)\leq x}}\la_S(\CO)h(\CO).
$$
Then as $x\ra\infty$ we have
$$
\pi_S(x)\sim\frac{e^{3x}}{3x}.
$$
More sharply,
$$
\pi_S(x)=L(3x)+O\left(\frac{e^{9x/4}}{x}\right)
$$
as $x\ra\infty$.

This result was again proved by means of a Prime Geodesic Theorem, this time for symmetric spaces formed as a compact quotient of the group $\SL_3(\R)$.  There exists a Selberg trace formula for such spaces (\cite{Wallach76}) by means of which the required properties of a generalised Selberg zeta function can be deduced (\cite{Deitmar00}) in order to prove the Prime Geodesic Theorem.  The class number formula is then deduced by means of a correspondence between primitive closed geodesics and orders in complex cubic number fields, under which the lengths of the geodesics correspond to the regulators of the number fields.

In this paper we follow the methods of \cite{Deitmar02}.
We make use of the prime geodesic theorem in \cite{primgeoSL4}.  Let us quote this result as follows.

Let $G=\SL_4(\R)$ and let $K$ be the maximal compact subgroup $\SO(4)$.  Let $\Ga\subset G$ be discrete and cocompact.  We then have a one to one correspondence between conjugacy classes in $\Ga$ and free homotopy classes of closed geodesics on the symmetric space $X_{\Ga}=\Ga\bs G/K$.  Let
$$
A^-=\left\{ \matrixfour{a}{a}{a^{-1}}{a^{-1}}:0<a<1\right\}
$$
and
$$
B=\matrixtwo {\SO(2)}{\SO(2)}.
$$
Let $\CE(\Ga)$ be the set of primitive conjugacy classes $[\ga]$ in $\Ga$ such that $\ga$ is conjugate in $G$ to an element $a_{\ga}b_{\ga}$ of $A^-B$.  For $\ga\in\Ga$ we write $a_{\ga}$ also for the top left entry in the matrix $a_{\ga}$ and define the length $l_{\ga}$ of $\ga$ to be $8\log a_{\ga}$.  Let $G_{\ga},\Ga_{\ga}$ be the centralisers of $\ga$ in $G$ and $\Ga$ respectively and let $K_{\ga}=K\cap G_{\ga}$.  For $x>0$ define the function
$$
\pi(x)=\sum_{{[\ga]\in\CE(\Ga)}\atop{e^{l_{\ga}}\leq x}}\chi_1(\Ga_{\ga}),
$$
where $\chi_1(\Ga_{\ga})$ is the first higher Euler characteristic of the symmetric space $X_{\Ga_{\ga}}=\Ga_{\ga}\bs G_{\ga}/K_{\ga}$.

\begin{theorem}\textnormal{(Prime Geodesic Theorem)}
For $x\rightarrow\infty$ we have
$$
\pi(x)\sim\frac{2x}{\log x}.
$$
More sharply,
$$
\pi(x)=2\,\li(x)+O\left(\frac{x^{3/4}}{\log x}\right)
$$
as $x\rightarrow\infty$, where $\li(x)=\int_2^x \rez{\log t} dt$ is the integral logarithm.
\end{theorem}

We now state our result on class numbers.  Let $S$ be a finite, non-empty set of prime numbers containing an even number of elements.  We define the sets $C(S)$ and $O(S)$ and the constants $\la_S(\CO)$ as above, replacing complex cubic fields with totally complex quartic fields in the definitions.  A totally complex quartic field has at most one real quadratic subfield, as can be seen by comparing numbers of fundamental units.  Let $C^c(S)\subset C(S)$ be the subset of fields with no real quadratic subfield and let $O^c(S)\subset O(S)$ be the subset of isomorphy classes of orders in fields in $C^c(S)$.  Let $R(\CO)$ and $h(\CO)$ once again denote respectively the regulator and class number.

\begin{theorem}\textnormal{(Main Theorem)}
\label{thm:main1}
For $x>0$ we define
$$
\pi_S(x)=\sum_{{\CO\in O^c(S)}\atop{R(\CO)\leq x}}\la_S(\CO)h(\CO).
$$
Then as $x\ra\infty$ we have
$$
\pi_S(x)\sim\frac{e^{4x}}{8x}.
$$
\end{theorem}
Note the surprising factor $1/2$!
We give an explanation below.

Transfering the methods of \cite{Deitmar02} to the case of totally complex quartic fields does not proceed entirely smoothly, there are various technical difficulties to be overcome.  For instance, the correspondence between primitive closed geodesics and orders is considerably more complicated than in either the real quadratic or the complex cubic case.  As a result of this extra complexity our result is weaker than that obtained in the complex cubic case.  We have had to impose the extra condition on the finite set $S$ of prime numbers that it must contain an even number of elements and we have been unable to provide an error term in our final asymptotic.

Also we have had to restrict ourselves to counting the orders in fields without a real quadratic subfield, as the fields with a real quadratic subfield cannot be counted using our method.  Indeed, any real quadratic field which occurs as a subfield of a totally complex quartic field may in fact occur as a subfield of infinitely many such fields.  The fundamental units in the quadratic field are powers of fundamental units in the quartic fields, so the (possibly infinitely many) quartic fields all have regulator less than or equal to that of the quadratic field.  Hence an asymptotic formula for a sum of class numbers of orders bounded by the regulators is not even possible in the case of totally complex quartic fields with a real subfield.  

In comparison with the real quadratic and complex cubic cases one might have expected to get the asymptotic
$
\pi_S(x)\sim\frac{e^{4x}}{4x}.
$
In the result that we do in fact get there is an extra factor of one half.  The fact that we have had to restrict the fields over which we are counting gives an explanation for this discrepancy.

We also prove the following asymptotic in which we introduce an extra factor into the summands.  For an order $\CO\in O(S)$ this extra factor is defined in terms of the arguments of the fundamental units of $\CO$ under the embeddings of $F$ into $\C$ as follows.  If $\ep$ is a fundamental unit in $\CO$ then $\ep^{-1},\ze\ep$ and $\ze\ep^{-1}$ are also fundamental units in $\CO$, where $\ze$ is a root of unity contained in $\CO$.  Let
$$
\nu(\CO)=\rez{2\mu_{\CO}}\sum_{\ep}\prod_{\al}\left(1-\frac{\al(\ep)}{|\al(\ep)|}\right),\index{$\nu(\CO)$}
$$
where $\mu_{\CO}$ is the number of roots of unity in $\CO$, the sum is over the $2\mu_{\CO}$ different fundamental units of $\CO$ and the product is over the embeddings of $\CO$ into $\C$.  We have:

\begin{theorem}
\label{thm:main2}
For $x>0$ we define
$$
\tilde{\pi}_S(x)=\sum_{{\CO\in O^c(S)}\atop{R(\CO)\leq x}}\nu(\CO)\la_S(\CO)h(\CO).
$$
Then as $x\ra\infty$ we have
$$
\tilde{\pi}_S(x)\sim\frac{e^{4x}}{2x}.
$$
\end{theorem}

An element $\ga$ in $\Ga$ is called \emph{regular} if it centraliser in $G$ is a torus and \emph{non-regular} otherwise.  The factor $\nu(\CO)$ was originally introduced in order to separate the contribution of non-regular elements from that of the regular elements in the Prime Geodesic Theorem.  As it turned out, the complexity of the correspondence between geodesics and orders meant that we did not really gain anything from this approach.  However, the factor $\nu(\CO)$ contains information about the arguments of the fundamental units in the order $\CO$ which is interesting in its own right.  In comparison with Theorem \ref{thm:main1} we can see that ``on average" the value of $\nu(\CO)$ as $R(\CO)$ goes to infinity is 4.  If $\ep$ is a fundamental unit in $\CO$ with arguments $\th,-\th,\phi,-\phi$ under the four embeddings of $\CO$ into $\C$ then a simple calculation shows that
$$
\prod_{\al}\left(1-\frac{\al(\ep)}{|\al(\ep)|}\right)=4(1-\cos 2\th)(1-\cos 2\phi).
$$
Since this takes ``on average" the value 4 we can see that there is a sense in which we can say that the arguments of the fundamental units in the orders $\CO$ are evenly distributed about $\pm\pi/2$ as $R(\CO)$ goes to infinity.

In what follows we give a summary of the arguments and techniques used in the proof of our main results.  In particular we shall point out the difficulties that arise in applying the methods of \cite{Deitmar02} to our situation.  

In Section \ref{ch:DivAlg} we use a classification of the set of equivalence classes of division algebras over $\Q$ by means of a description of the Brauer group of $\Q$ (\cite{Pierce82}, Theorem 18.5) to show the existence of a division algebra $M(\Q)$ whose maximal subfields include the fields in the set $C^c(S)$ and such that $M(\Q)\ox\R\cong\Mat_4(\R)$.  In fact the maximal subfields of $M(\Q)$ are all quartic extensions of $\Q$ and the set of totally complex maximal subfields of $M(\Q)$ coincides with $C(S)$.  We further show that for an order $\CO$ in $O(S)$ the number of embeddings of $\CO$ into $M(\Q)$, up to conjugation by the unit group of the maximal order $M(\Z)$, is $\la_S(\CO)h(\CO)$.  The restriction on the set $S$ of prime numbers that it has to contain an even number of elements is a consequence of the classification of division algebras over $\Q$.

In Section \ref{ch:Compare} we prove a correspondence between the primitive, closed geodesics in the Prime Geodesic Theorem and the orders of totally complex quartic fields.  Under this correspondence the lengths of the geodesics correspond to the regulators of the orders.  It turns out that it is the regular, weakly neat elements of $\CE(\Ga)$ which correspond to orders in totally complex quartic fields with no real quadratic subfield.  We use a ``sandwiching" argument to isolate these elements in the Prime Geodesic Theorem.

We finish this introduction with a few remarks about the limitations of the method used here in terms of further applications and mention a couple of other recent results in the same direction.  In order to be able to make use of the trace formula for compact spaces we have had to limit our sum over class numbers by means of the choice of a finite set of primes, as described above.  In \cite{DeitmarHoffmann05} Deitmar and Hoffmann have been able to use a different trace formula to prove that as $x\ra\infty$
$$
\sum_{R(\CO)\leq x}h(\CO)\sim\frac{e^{3x}}{3x},
$$
where the sum is over all isomorphism classes of orders in complex cubic fields.

In order to get the error term in the Prime Geodesic Theorem we have made use of the classification of the unitary dual of $\SL_4(\R)$.  At present the unitary dual is not known for any higher dimensional groups so a Prime Geodesic Theorem with error term is not possible using our methods.  Finally we mention that the correspondence between geodesics and orders actually works by identifying primitive geodesics with fundamental units in orders.  By Dirichlet's unit theorem, an order in a number field $F$ has a unique fundamental unit (up to inversion and multiplication by a root of unity) only if $F$ is real quadratic, complex cubic or totally complex quartic.  Hence an asymptotic of our form can be proven only in these three cases.  In \cite{Deitmar04}, the first named author has proved a Prime Geodesic Theorem for higher rank spaces, from which he deduces an asymptotic formula for class numbers of orders in totally real number fields of prime degree.

\newpage
\section{Division Algebras of Degree Four}
    \label{ch:DivAlg}
\subsection{Central simple algebras and orders}
Let $F$ be a field.    Let $\CCC(F)$ \index{$\CCC(F)$} denote the set of isomorphism classes of finite dimensional central simple algebras over $F$ and $\CCD(F)$ \index{$\CCD(F)$} the set of isomorphism classes of finite dimensional division algebras with centre $F$, then $\CCD(F)\subset\CCC(F)$.

For convenience we collect together here a number of facts about central simple algebras and their orders, which will be used in the sequel.  Let $F$ be a field, $K$ an extension of $F$ and $A$ a finite dimensional algebra over $F$.

\begin{proposition}
\label{pro:algebras}

(a) If $n=\dim_F A$, then there exists an injective homomorphism $A\hra\Mat_n(F)$.

(b) $\dim_K(A\ox K)=\dim_F A$.

(c) If $A\in\CCC(F)$, then $A\ox K\in\CCC(K)$.

(d) If $B$ is a simple subalgebra of $A\in\CCC(F)$ with centraliser $C_A(B)$ in A, then $$(\dim_F B)(\dim_F C_A(B))=\dim_F A.$$
\end{proposition}
\prf
All the statements of the proposition are proved in \cite{Pierce82}.  Statement (a) is Corollary 5.5b; (b) is Lemma 9.4; (c) is Proposition 12.4b(ii) and (d) is Theorem 12.7.
\qed

The next proposition gives some results about subfields of algebras.

\begin{proposition}
\label{pro:subfields}
(a) If $A\in\CCD(F)$, then for every $x\in A$ the set $F[x]=\left\{f(x):f\in F[X]\right\}$ is a subfield of $A$.

(b) If $A\in\CCC(F)$, then $\dim_F A =m^2$ for some $m\in\N$.  The natural number $m$ is called the degree of $A$.  If $E$ is a subfield of $A$ then $[E:F]$ divides $m$.

If $K$ is a subfield of $A$ containing $F$ then it is said to be \emph{strictly maximal} \index{strictly maximal subfield} if $[K:F]=\deg A$.

(c) A subfield $K$ of $A\in\CCC(F)$ is strictly maximal if and only if $C_A(K)=K$.  If $A\in\CCD(F)$ then every maximal subfield of $A$ is strictly maximal.

(d) If $A\in\CCC(F)$ and $[K:F]=\deg A=n$, then $A\ox K\cong\Mat_n(K)$ if and only if $K$ is isomorphic as an $F$-algebra to a strictly maximal subfield of $A$.
\end{proposition}
\prf
All the statements of the proposition are proved in \cite{Pierce82}.  Statement (a) is Lemma 13.1b; (b) is Corollary 13.1a; (c) is Corollary 13.1b and (d) is Corollary 13.3.
\qed

Now let $F$ be a number field and $\CO_F$ the ring of integers in $F$.  For any finite dimensional vector space $V$ over $F$, a \emph{full $\CO_F$-lattice} \index{full $\CO_F$-lattice} in $V$ is a finitely generated $\CO_F$-submodule $M$ in $V$ which contains a basis of $V$.  A subring $\CO$ of $A$ which is also a full $\CO_F$-lattice in $A$ is called an \emph{$\CO_F$-order}\index{order!$\CO_F$-order}, or simply an \emph{order}\index{order}, in $A$.  A \emph{maximal} \index{order!maximal} $\CO_F$-order in $A$ is an $\CO_F$-order which is not properly contained in any other $\CO_F$-order in $A$.  An element $a\in A$ is said to be \emph{integral} \index{integral element} over $\CO_F$ if it is a root of a monic polynomial with coefficients in $\CO_F$.  The \emph{integral closure} \index{integral closure} of $\CO_F$ in $A$ is the set of all elements of $A$ which are integral over $\CO_F$.

\begin{proposition}
\label{pro:orders}
(a)(Skolem-Noether Theorem) \index{Skolem-Noether Theorem} Let $A\in\CCC(F)$, and $B$ a simple subring of $A$ such that $F\subset B\subset A$.  Then every $F$-isomorphism of $B$ onto a subalgebra of $A$ extends to an inner automorphism of $A$.

(b) Every element of an $\CO_F$-order $\CO$ is integral in $A$ over $\CO_F$.

(c) The $\CO_F$-order $\CO$ is maximal in $A$ if and only if for each prime ideal $\p$ of $\CO_F$ the $\p$-adic completion $\CO_{\p}$ is a maximal $\CO_{F,\p}$-order in $A_{\p}$.

(d) If $A\in\CCD(F_{\p})$ for some prime ideal $\p$ of $\CO_F$, then the integral closure of $\CO_{F,\p}$ in $A$ is the unique maximal $\CO_{F,\p}$-order in $A$.
\end{proposition}
\prf
All the statements of the proposition are proved in \cite{Reiner75}.  Statement (a) is Theorem 7.21; (b) is Theorem 8.6; (c) is Corollary 11.6 and (d) is Theorem 12.8.
\qed

\subsection{Division algebras of degree four}

Let $A$ be a central simple algebra over a field $F$.  By the Wedderburn Structure Theorem \index{Wedderburn Structure Theorem} (\cite{Pierce82}, Theorem 3.5) the algebra $A$ is isomorphic to $\Mat_n(D)$, where $n\in\N$ and $D$ is a finite dimensional division algebra over $F$, which by \cite{Pierce82}, Proposition 12.5b is unique up to isomorphism.  If $B$ is another central simple algebra over $F$ isomorphic to $\Mat_m(E)$, where $m\in\N$ and $E$ is a finite dimensional division algebra over $F$, then $A$ and $B$ are called \emph{Morita equivalent} \index{Morita equivalent} if and only if $D\cong E$.  By \cite{Pierce82}, Proposition 12.5a, the tensor product over $F$ induces a group structure on the set of equivalence classes.  The group defined in this way is called the \emph{Brauer group}.\index{Brauer group}

By \cite{Pierce82}, Theorem 17.10, if $F$ is a local, non-archimedean field, then there is a canonical isomorphism between the Brauer group of $F$ and the group $\Q/\Z$.  The \emph{Brauer invariant} \index{Brauer invariant} $\Inv A$ \index{$\Inv A$} of $A$ is defined to be the image in $\Q/\Z$ of the Morita equivalence class of $A$ under this isomorphism.  The \emph{Schur index} \index{Schur index} $\Ind A$ \index{$\Ind A$} of $A$ is defined to be the degree of $D$.  By \cite{Pierce82}, Corollary 17.10a(iii), the order of the Brauer invariant of $A$ in $\Q/\Z$ equals $\Ind A$.

The only (associative) division algebras over $\R$ are $\R$ itself, $\C$ and the quaternions $\H$ (see \cite{Pierce82}, Corollary 13.1c).  Since $\C$ is itself a field, any algebra which is central simple over $\R$ cannot be isomorphic to a matrix algebra over $\C$.  The Brauer group of $\R$ is therefore isomorphic to $\Z/2\Z$.  Let $A$ be a central simple algebra over $\R$.  If $A\cong\Mat_n(\R)$ for some $n\in\N$, then the Brauer invariant of $A$ is defined to be zero and the Schur index of $A$ is defined to be one.  If $A\cong\Mat_n(\H)$ for some $n\in\N$, then the Brauer invariant of $A$ is one half and the Schur index of $A$ is two.

Let $M$ be a division algebra of degree 4 over $\Q$.  Fix a maximal order $M(\Z)$ \index{$M(\Z)$} in $M$.  If $R$ is a commutative ring with unit then we denote by $M(R)$ \index{$M(R)$} the $R$-algebra $M(\Z)\ox_{\Z} R$.  In particular we have $M\cong M(\Q)$ \index{$M(\Q)$}.  For any ring $R$ the reduced norm induces a map $\det\!\!:M(R)\ra R$.  (We note that in the case that $M(R)\cong\Mat_4(R)$ the reduced norm of an element of $M(R)$ equals its determinant, justifying the choice of notation; see \cite{Pierce82}, Chapter 16.)

Let $p$ be a prime.  From Proposition \ref{pro:algebras}(b) we know that $\deg M(\Q_p)=\deg M(\Q)=4$.  It then follows that if $\Inv M(\Q_p)$ is equal to $\rez{4}$ or $\frac{3}{4}$ then $M(\Q_p)$ is a division algebra of degree four over $\Q_p$, and if $\Inv M(\Q_p)=0$ then $M(\Q_p)\cong\Mat_4(\Q_p)$.  The other possibility, that $\Inv M(\Q_p)=\rez{2}$, does not interest us here.  Proposition \ref{pro:algebras}(b) also tells us that $\deg M(\R)=4$.  We know that $\Inv M(\R)$ equals zero or one half.  It follows that $M(\R)$ is isomorphic to $\Mat_4(\R)$ or $\Mat_2(\H)$ respectively.

Let $S$ \index{$S$} be a finite, non-empty set of prime numbers with an even number of elements.  We say that $M(\Q)$ \emph{splits} \index{splitting!of division algebra} over a prime $p$ if $M(\Q_p)\cong\Mat_4(\Q_p)$.  For all primes $p$ we define $i_p$ as follows.  If $p\in S$ define $i_p$ to be either $\rez{4}$ or $\frac{3}{4}$ in such a way that $\sum_{p\in S}i_p\in\Z$.  Note that we are able to do this since we have specified that $S$ must contain an even number of elements.  For all other $p$ define $i_p=0$.  Then \cite{Pierce82}, Theorem 18.5 tells us that we may choose $M$ so that $\Inv M(\Q_p)=i_p$ for all primes $p$ and $\Inv M(\R)=0$.  We can then see that the set of places at which $M(\Q)$ does not split coincides with the set $S$.  More particularly, for $p\in S$ we have that $M(\Q_p)$ is a division $\Q_p$-algebra, for $p\notin S$ we have that $M(\Q_p)\cong\Mat_4(\Q_p)$ and we also have $M(\R)\cong\Mat_4(\R)$.

For a commutative ring $R$ with unity, let \index{$\CG(R)$}
$$
\CG(R)=\left\{x\in M(R):\det(x)=1\right\}.
$$
Then $\CG$ \index{$\CG$} is a linear algebraic group, defined over $\Z$.  Let $\Ga=\CG(\Z)$, \index{$\Ga$} then $\Ga$ forms a discrete subgroup of $G=\CG(\R)\cong\SL_4(\R)$. Since $M(\Q)$ is a division algebra, it follows that $\CG$ is anisotropic over $\Q$ and so $\Ga$ is cocompact in $G$ (see \cite{BorelHarder78}, Theorem A).  Let $P$ be the parabolic subgroup \index{$P$}
$$
P=\matrix{*}{*}{\begin{array}{cc} 0 & 0 \\ 0 & 0 \end{array}}{*}
$$
of G.  Then $P=MAN$, where \index{$M$}
\begin{eqnarray*}
M & = & {\rm S}\matrixtwo{\SL_2^{\pm}(\R)}{\SL_2^{\pm}(\R)} \\
  & = & \left\{\matrixtwo{X}{Y}:X,Y\in\Mat_2(\R),\ \det X=\det Y=\pm 1\right\},
\end{eqnarray*}
$A=\left\{ \diag (a,a,a^{-1},a^{-1}):a>0\right\}$ \index{$A$} and the elements of $N$ \index{$N$} have ones on the diagonal and the only other non-zero entries in the top right two by two square.  Let \index{$B$}
$$
B=\matrixtwo{\SO(2)}{\SO(2)}.
$$
Then $B$ is a compact subgroup of $M$.

Let $A^-=\left\{ \diag(a,a,a^{-1},a^{-1}):a\in (0,1)\right\}$ \index{$A^-$} be the negative Weyl chamber of $A$.  Let $\CE_P(\Ga)$ \index{$\CE_P(\Ga)$} be the set of conjugacy classes $[\ga]$ in $\Ga$ such that $\ga$ is conjugate in $G$ to an element $a_{\ga}b_{\ga}$ of $A^- B$ and let $\CE_P^p(\Ga)$ \index{$\CE_P^p(\Ga)$} be the subset of primitive conjugacy classes.

We say $g\in G$ is \emph{regular} \index{regular} if its centraliser is a torus and \emph{non-regular} \index{non-regular} (or \emph{singular}) \index{singular} otherwise.  Clearly, for $\ga\in\Ga$ regularity is a property of the $\Ga$-conjugacy class $[\ga]$, we denote by $\CE_P^{p,\reg}(\Ga)$ \index{$\CE_P^{p,\reg}(\Ga)$} the regular elements of $\CE_P^p(\Ga)$.  By abuse of notation we sometimes write $\ga\in\CE_P^p(\Ga)$ or $\ga\in\CE_P^{p,\reg}(\Ga)$ when we mean the conjugacy class of $\ga$.

We will call a quartic field extension $F/\Q$ \emph{totally complex} \index{totally complex field} if it has two pairs of conjugate complex embeddings and no real embeddings into $\C$.

\subsection{Subfields of $M(\Q)$ generated by $\Ga$}
\begin{lemma}
\label{lem:centralisers}
Let $[\ga]\in\CE_P(\Ga)$.

The centraliser $M(\Q)_{\ga}$ of $\ga$ in $M(\Q)$ is a totally complex quartic field if and only if $\ga$ is regular.

The centraliser $M(\Q)_{\ga}$ is a quaternion algebra over the real quadratic field $\Q[\ga]$ if and only if $\ga$ is non-regular.
\end{lemma}

\prf
The centraliser $M(\Q)_{\ga}$ of $\ga$ in $M(\Q)$ is a division subalgebra of $M(\Q)$.  Suppose that $M(\Q)_{\ga}=\Q$.  Then since $\ga\in M(\Q)_{\ga}$ we must have $\ga\in\Q$, but then $\ga$ is central so $M(\Q)_{\ga}=M(\Q)$, a contradiction.  Suppose instead that $M(\Q)_{\ga}=M(\Q)$, ie. $\ga$ is central.  Then $\ga\in\Q$, so $\det\ga=\ga^4=1$ and it follows that $\ga=\pm 1$, which possibility is excluded since $\ga=\pm 1$ is not primitive.

By Proposition \ref{pro:algebras}(d) we have that $\dim_{\Q}M(\Q)_{\ga}$ divides $\dim_{\Q}M(\Q)$, which is 16, so $\dim_{\Q}M(\Q)_{\ga}$ is equal to 2, 4 or 8.  Also
$$
\dim_{\Q}M(\Q)_{\ga}=\dim_{\R}M(\R)_{\ga}=\dim_{\R}M(\R)_{a_{\ga}b_{\ga}}.
$$
Depending on whether neither, one or both of the $\SO(2)$ components of $b_{\ga}$ are $\pm I_2$, the dimension $\dim_{\R}M(\R)_{a_{\ga}b_{\ga}}$ is equal to 4, 6 or 8 respectively.  Hence $\dim_{\Q}M(\Q)_{\ga}$ is either 4 or 8.  If $\dim_{\Q}M(\Q)_{\ga}=4$ then $b_{\ga}=(R(\th),R(\phi))$ with $\th,\phi\notin\pi\Z$, where
$$
R(\th)=\matrix{\cos\th}{-\sin\th}{\sin\th}{\cos\th}.
$$
If $\dim_{\Q}M(\Q)_{\ga}=8$ then $b_{\ga}$ is diagonal.

In the first case $G_{\ga}\cong AB\cong\R^+\x\SO(2)\x\SO(2)$ so $\ga$ is regular.  
From Proposition \ref{pro:algebras}(d) and Proposition \ref{pro:subfields}(a) it follows that $\Q[\ga]$ is a subfield of $M(\Q)$ of degree 4.  
Let $f_{\ga}(x)$ be the minimal polynomial of $\ga$ over $\Q$.  
Then $a_{\ga}b_{\ga}$ also satisfies $f_{\ga}(x)=0$.  Now $a_{\ga}=\diag(a,a,a^{-1},a^{-1})$ for some $a\in(0,1)$ and $b_{\ga}=(R(\th),R(\phi))$ for some $\th, \phi\notin\pi\Z$, so the complex numbers $z_1=ae^{i\th}$ and $z_2=a^{-1}e^{i\phi}$ also satisfy $f_{\ga}(x)=0$.
Neither $z_1$ nor $z_2$ are real, nor do we have $z_1=\bar{z_2}$, so $\Q[\ga]\cong\Q[z_1,z_2]$, which is a totally complex field.  
The division subalgebra $M(\Q)_{\ga}$ of $M(\Q)$ has dimension four over $\Q$ and contains $\Q[\ga]$, hence is equal to $\Q[\ga]$.

Suppose instead that $b_{\ga}$ is diagonal.  Then $G_{\ga}\cong AM$, so $\ga$ is not regular.  From Proposition \ref{pro:algebras}(d) and Proposition \ref{pro:subfields} (a) it follows that $\Q[\ga]=\left\{f(\ga):f(x)\in \Q[x]\right\}$ is a subfield of $M(\Q)$ of degree 2.  Let $f_{\ga}(x)$ be the minimal polynomial of $\ga$ over $\Q$.  Then $a_{\ga}b_{\ga}$ also satisfies $f_{\ga}(x)=0$.  Now $a_{\ga}=\diag(a,a,a^{-1},a^{-1})$ for some $a\in(0,1)$ and $b_{\ga}=(\pm I_2,\pm I_2)$.  Let us suppose that $b_{\ga}=I_4$, the other cases are similar.  Then $f_{\ga}(a)=f_{\ga}(a^{-1})=0$, so $f_{\ga}=(x-a)(x-a^{-1})$ and $\Q[\ga]$ is a real quadratic field.

The division algebra $M(\Q)_{\ga}$ is central simple over its centre $Z(M(\Q)_{\ga})$.  Clearly, $\Q[\ga]$ is contained in $Z(M(\Q)_{\ga})$.  By Proposition \ref{pro:subfields}(b), the dimension of $M(\Q)_{\ga}$ over its centre is $m^2$ for some $m\in\N$.  Hence either $M(\Q)_{\ga}$ is a field or $Z(M(\Q)_{\ga})=\Q[\ga]$.  However, by Proposition \ref{pro:subfields}(b) again, if $F$ is a subfield of $M(\Q)$ then $[F:\Q]$ = 1, 2 or 4.  This rules out the possibility of $M(\Q)_{\ga}$ being a field, since $\dim_{\Q}M(\Q)_{\ga}=8$.  We conclude that $Z(M(\Q)_{\ga})=\Q[\ga]$ and by \cite{Pierce82}, Theorem 13.1, $M(\Q)_{\ga}$ is a quaternion algebra over $\Q[\ga]$.
\qed

\subsection{Field and order embeddings}

In this section we prove a number of lemmas which will be needed later.

\begin{lemma}
\label{lem:unit}
Let $u\in M(\Z)$.  Then $u$ is a unit in $M(\Z)$ iff $\det(u)=\pm 1$.  If $\Q[u]$ is quadratic or is totally complex quartic, or if $u=\pm 1$, then $\det(u)=1$.
\end{lemma}

\prf
By Proposition \ref{pro:subfields}(a) the set $F=\Q[u]$ is a subfield of $M(\Q)$ containing $u$.  The set $F\cap M(\Z)$ is an order in $F$.  The element $u$ is in $F\cap M(\Z)$ so by Proposition \ref{pro:orders}(b) it is in the integral closure $\CO_F$ of $\Z$ in $F$.  

The field norm $N_{F|\Q}:F\ra\Q$ can be written as a product over all embeddings of $F$ into the complex field,
$
N_{F|\Q}(x)=\prod_{\si}\si x,
$
and it is standard that $u$ is a unit in $\CO_F$ if and only if $N_{F|\Q}(u)=\pm 1$.  Proposition 16.2a of \cite{Pierce82} tells us that
\begin{equation}
\label{eqn:detNormReln}
\det(u)=N_{F|\Q}(u)^k,
\end{equation}
where $\deg M(\Q)=k[F:\Q]$.  The first statement of the lemma follows.

If $u=\pm 1$ (so that $F=\Q$) or $F$ is quadratic then by (\ref{eqn:detNormReln}) we have $\det(u)=1$.  Suppose $F$ is a totally complex quartic extension.  Let $f(X)=X^4+a_3 X^3+a_2 X^2+a_1 X+a_0$ be the minimal polynomial of $u$ over $\Q$.  Then $\det(u)=a_0$.  If $a_0=\det(u)=-1$ then $f$ must have a real root, which contradicts the supposition on $F$.  Hence $\det(u)=1$.
\qed

The following lemma holds for any algebraic extension $F$ of $\Q$.
\begin{lemma}
\label{lem:maxorders}
Let $\CO$ be an order in the number field $F/\Q$ and $p$ a prime number.  Let $F_p=F\ox\Q_p$ \index{$F_p$} and $\CO_p=\CO\ox\Z_p$\index{$\CO_p$}.  Then $\CO_p$ is the maximal order of $F_p$ for all but finitely many primes $p$.
\end{lemma}

\prf
Let $\CO_M$ be the maximal order of $F$ and let $\CO_{M,p}=\CO_M\ox\Z_p$.  We define the \emph{conductor} \index{conductor} of $\CO$ to be the ideal $\CF=\{\al\in\CO_M|\al\CO_M\subset\CO\}$ of $\CO_M$.  Then for all primes $p$, if $\CF\nsubseteq p\CO_M$, then $\CO_p=\CO_{M,p}$.  In fact, if $\CF\nsubseteq p\CO_M$ then there is an element $\al\in\CF$ such that $\al\notin p\CO_M$.  Every element of $\CO_{M,p}$ can be written as a unit of $\CO_{M,p}$ times a power of $p$.  Since $\al\notin p\CO_M$ we have $\al\in\CO_{M,p}^{\x}$, which implies $\al\CO_{M,p}=\CO_{M,p}$.  Then, as $\al\in\CF$ it follows that
$$
\CO_{M,p}=\al\CO_{M,p}\subset\CO_p.
$$
It is clear that $\CO_p\subset\CO_{M,p}$.  Since $\CF\subset p\CO_M$ for only finitely many primes $p$, the lemma follows.
\qed

Let $F$ be a number field.  A prime number $p$ is called {\it non-decomposed} \index{prime!non-decomposed} \index{non-decomposed (prime)} in $F$ if there is only one place in $F$ lying above $p$.

\begin{lemma}
\label{lem:fieldembedding}
A number field $F$ embeds into $M(\Q)$ if and only if $[F:\Q]=1, 2$ or $4$ and $p$ is non-decomposed in $F$ for all $p\in S$.
\end{lemma}

\prf
The statement about the degree of $F$ is Proposition \ref{pro:algebras} (c).  The case $F=\Q$ is trivial.

Let $F$ be a quartic field and suppose that $F$ embeds into $M(\Q)$, then $F$ is a strictly maximal subfield of $M(\Q)$.  We say that a subfield $K$ of $M(\Q)$ is a \emph{splitting field} \index{splitting!field} for $M(\Q)$ (or that $K$ \emph{splits} $M(\Q)$) if $M(\Q)\ox K\cong \Mat_4(K)$.  Proposition \ref{pro:subfields}(d) tells us that a quartic field embeds into $M(\Q)$ if and only if it splits $M(\Q)$.

For $p$ a rational prime, the field $\Q_p$ is an extension of $\Q$.  Hence, by Proposition \ref{pro:algebras}(c), the $\Q_p$-algebra $M(\Q_p)$ is central simple over $\Q_p$.  Let $p\in S$.  The Schur index, $\Ind M(\Q_p)$, of $M(\Q_p)$ is, by \cite{Pierce82}, Corollary 17.10a(iii), equal to the order of $[M(\Q_p)]$ in the Brauer group of $\Q_p$, where the square brackets denote the Morita equivalence class.  By choice of $M(\Q)$ we have $\Ind M(\Q_p)=4$.  Let $\p_i$ be the (finitely many) primes in $F$ above $p$.  Then Theorem 32.15 of \cite{Reiner75} tells us that, since $F$ splits $M(\Q)$,
$$
4=\Ind M(\Q_p)\ |\ [F_{\p_i}:\Q_p]
$$
for all $i$.  Corollary 8.4 of Chapter II of \cite{Neukirch99} says that
\begin{equation}
\label{eqn:padicSum}
[F:\Q]=\sum_i[F_{\p_i}:\Q_p],
\end{equation}
which implies that there is only one prime $\p$ in $F$ above $p$ and that $[F_{\p}:\Q_p]=4$.

Conversely, suppose that for each $p\in S$ there is only one prime $\p$ in $F$ above $p$.  From (\ref{eqn:padicSum}) we get that $[F_{\p}:\Q_p]=4$, so
\begin{equation}
\label{eqn:indexDivides}
\Ind M(\Q_p)\ |\ [F_{\p}:\Q_p].
\end{equation}
For $p\notin S$, by choice of $M(\Q)$, the Schur index, $\ind M(\Q_p)$, is equal to 1, so in this case (\ref{eqn:indexDivides}) also holds for each $\p$ in $F$ above $p$.  Then, by \cite{Reiner75}, Theorem 32.15, the field $F$ splits $M(\Q)$.  By Proposition \ref{pro:subfields}(d) we then deduce that $F$ embeds into $M(\Q)$.

Now let $F$ be a quadratic field and suppose that $F$ embeds into $M(\Q)$, then we consider $F$ as a subfield of $M(\Q)$.  By Proposition \ref{pro:subfields}(c), the subfield $F$ of $M(\Q)$ is not maximal so there exists a maximal subfield $K$ of $M(\Q)$ properly containing $F$, which by the same proposition is quartic.  It was shown above that every prime in $S$ is non-decomposed in $K$ and hence also in $F$.

Suppose $p$ is non-decomposed in $F$ for all $p\in S$.  We choose a quaternion algebra $A$ over $F$ by specifying that the local Brauer invariant at $\p$ be $\rez{2}$ for all places $\p$ in $F$ over some $p\in S$, and that the local Brauer invariant be 0 at all other places of $F$.  Proposition \ref{pro:subfields}(c) tells us that there exists a subfield $K$ of $A$ containing $F$ with $[K:F]=2$.  Then by the same argument as above, if $\p$ is a place of $F$ over $p$ for some $p\in S$, then $\p$ is non-decomposed in $K$.  It follows that $p$ is non-decomposed in $K$ for all $p\in S$.  It was shown above that the quartic extension $K$ of $\Q$, and thus also the subfield $F$ of $K$, embeds into $M(\Q)$.  The lemma is proven.
\qed

Let $\A_{\fin}$ \index{$\A_{\fin}$} denote the ring of finite adeles over $\Q$, ie.
$$
\A_{\fin}=\left.\prod_p\right.'\Q_p,
$$
where $p$ ranges over the rational primes and $\prod'$ denotes the restricted product, that is, almost all components of a given element are integral.  Let $\hat{\Z}=\prod_p\Z_p\subset\A_{\fin}$\index{$\hat{\Z}$}.  There is a natural embedding of $\Q$ into $\A_{\fin}$, which maps $q\in\Q$ to the adele every one of whose components is $q$.  From this embedding we derive an embedding of $M(\Q)$ in $M(\A_{\fin})$ and an embedding of $\CG(\Q)$ in $\CG(\A_{\fin})$.

\begin{lemma}
\label{lem:Adeles}
$M(\A_{\fin})^{\x}=M(\hat{\Z})^{\x}\,M(\Q)^{\x}$.
\end{lemma}
\prf
We have the following commutative diagram with exact rows:
\begin{displaymath}
\xymatrix{1 \ar[r] & \CG(\Q) \ar[r] \ar[d] & M(\Q)^{\x} \ar[r]^{\det} \ar[d] & \Q^{\x} \ar[r] \ar[d] & 1 \\
   1 \ar[r] & \CG(\A_{\fin}) \ar[r] & M(\A_{\fin})^{\x} \ar[r]^{\det} & \A_{\fin}^{\x} \ar[r] & 1 \\
   1 \ar[r] & \CG(\hat{\Z}) \ar[r] \ar[u] & M(\hat{\Z})^{\x} \ar[r]^{\det} \ar[u] & \hat{\Z}^{\x} \ar[r] \ar[u] & 1}
\end{displaymath}
where the vertical arrows denote the natural embeddings.  The commutativity of the diagram and the exactness of the rows are clear, except that the surjectivity of the three determinant maps requires justification.  The surjectivity of the map $\det\!\!:M(\Q)^{\x}\ra\Q^{\x}$ follows from \cite{Reiner75}, Theorem 33.15.  For all primes $p$ we have that $M(\Q_p)$ is a central simple $\Q_p$-algebra so it follows from \cite{Reiner75}, Theorem 33.4 that the map $\det\!\!:M(\Q_p)^{\x}\ra\Q_p^{\x}$ is surjective.  For $p\in S$, we have that $M(\Q_p)$ is a division $\Q_p$-algebra and so an element $x$ of $M(\Q_p)$ is in $M(\Z_p)$ if and only if $\det x\in\Z_p$ (\cite{Reiner75}, Theorem 12.5).  For $p\notin S$, we have that $M(\Q_p)\cong\Mat_4(\Q_p)$ and
$$
\det\matrixfour{x}{1}{1}{1}=x
$$
for all $x\in\Z_p$.  It follows that the map $\det\!\!:M(\hat{\Z})^{\x}\ra\hat{Z}^{\x}$ is surjective.  From this we can then see that map $\det\!\!:M(\A_{\fin})^{\x}\ra\A_{\fin}^{\x}$ is also surjective.

Let $v\in M(\A_{\fin})^{\x}$ and let $w=\det v\in\A_{\fin}^{\x}$.  Since $\A_{\fin}^{\x}=\hat{\Z}^{\x}\Q^{\x}$, there exist $r\in\hat{\Z}^{\x}$ and $q\in\Q^{\x}$ such that $w=rq$.  By the surjectivity of $\det$, there exist $\bar{r}\in M(\hat{\Z})^{\x}$ and $\bar{q}\in M(\Q)^{\x}$ such that $\det\bar{r}=r$ and $\det\bar{q}=q$.  Let $\bar{v}=\bar{r}^{-1}v\bar{q}^{-1}$.  Then $\det\bar{v}=1$, so $\bar{v}\in\CG(\A_{\fin})$.  By the Strong Approximation Theorem (see \cite{Kneser66}), we have that $\CG(\Q)$ is dense in $\CG(\A_{\fin})$.  Hence there exist $\hat{r}\in\CG(\hat{\Z})$ and $\hat{q}\in\CG(\Q)$ such that $\bar{v}=\hat{r}\hat{q}$.  Finally, we have
$$
v=\bar{r}\bar{v}\bar{q}=\bar{r}(\hat{r}\hat{q})\bar{q}=(\bar{r}\hat{r})(\hat{q}\bar{q}),
$$
where $\bar{r}\hat{r}\in M(\hat{\Z})^{\x}$ and $\hat{q}\bar{q}\in M(\Q)^{\x}$.  The lemma follows.
\qed

\begin{lemma}
\label{lem:orderembedding}
Let $F$ be a field that embeds into $M(\Q)$.  For each embedding $\si:F\ra M(\Q)$ the order $\CO_{\si}=\si^{-1}\left(\si(F)\cap M(\Z)\right)$ \index{$\CO_{\si}$} is maximal at all $p\in S$.  Conversely, if $\CO$ is an order in the field $F$ which is maximal at all $p\in S$, then there exists an embedding $\si:F\ra M(\Q)$ such that $\CO=\CO_{\si}$.
\end{lemma}

Let $\si:F\ra M(\Q)$ be an embedding of the field $F$ into $M(\Q)$.  By Lemma \ref{lem:fieldembedding} each prime $p$ in the set $S$ is non-split in $F$.  Since for $p\in S$ there is only one place in $F$ above $p$, we can note that the $p$-adic completion $F_p=F\ox\Q_p$ is again a field (see \cite{Neukirch99}, Chapter 1, Proposition 8.3).  Firstly we need to show that for all $p\in S$, the completion $\CO_{\si,p}=\CO_{\si}\ox\Z_p$ is the maximal order in $F_p$.  By Proposition \ref{pro:orders}(c), $M(\Z_p)$ is the maximal $\Z_p$-order in $M(\Q_p)$.  Since $M(\Q_p)$ is a division algebra, $M(\Z_p)$ coincides with the integral closure of $\Z_p$ in $M(\Q_p)$, by Proposition \ref{pro:orders}(d).  We can extend $\si$ to an embedding of $F_p$ in $M(\Q_p)$ and then $\CO_{\si,p}=\si^{-1}\left(\si(F_p)\cap M(\Z_p)\right)$ is the integral closure of $\Z_p$ in $F_p$, which is the maximal order of $F_p$.

For the converse, let $\CO$ be an order of $F$ such that the completion $\CO_{p}$ is maximal for each $p\in S$.  Via $\si$ we consider $F$ to be a subfield of $M(\Q)$.  For any $u\in M(\Q)^{\x}$ let $\CO_u =F\cap u^{-1}M(\Z)u$.  Let $^u\si=u\si u^{-1}$, that is: $^u\si(x)=u\si(x)u^{-1}$ for all $x\in F$.  We will show that there is a $u\in M(\Q)^{\x}$ such that $\CO=\CO_u$.  The embedding $^u\si$ is then the one required by the lemma.

Let $\CO_1=F\cap M(\Z)$.  Since $\CO$ and $\CO_1$ are orders, by Lemma \ref{lem:maxorders} they are both maximal at all but finitely many places.  Let $T$ be the set of primes $p$ such that either $\CO$ or $\CO_1$ is not maximal at $p$.  Then $T$ is finite and $T\cap S=\varnothing$.  Furthermore, if $T_F$ denotes the set of places of $F$ lying over $T$, we have that for any place $\p$ of $F$ with $\p\notin T_F$ the completions $\CO_{\p}$ and $\CO_{1,\p}$ are maximal in $F_{\p}$ and thus equal, by uniqueness of maximal orders.  Hence, for $p\notin T$ we have
$$
\CO_p=F_p\cap M(\Z_p).
$$

Let $p\in T$, then $p\notin S$ so $M(\Q_p)\cong\Mat_4(\Q_p)$.  Considering $F_p$ as a $\Q_p$-vector space we see that we can embed $F_p$ into $\Lin_{\Q_p}(F_p,F_p)$ by sending an element $x\in F_p$ to the linear map $a\mapsto ax$.  We choose a $\Z_p$-basis of $\CO_p$, which is then also a $\Q_p$-basis of $F_p$.  This gives us an isomorphism
$$
\Lin_{\Q_p}(F_p,F_p)\cong\Mat_4(\Q_p)\cong M(\Q_p).
$$
This isomorphism and the above embedding $F_p\hra\Lin_{\Q_p}(F_p,F_p)$ then give us an embedding $\si_p:F_p\hra M(\Q_p)$ such that $\CO_p=\si_p^{-1}\left(\si_p(F_p)\cap M(\Z_p)\right)$.  The map $\si_p$ is a $\Q_p$-isomorphism of $F_p$ (considered as a subfield of $M(\Q_p)$ via a suitable extension of $\si$) onto its image in $M(\Q_p)$, so by the Skolem-Noether Theorem (Proposition \ref{pro:orders}(a)), there exists $u_p\in M(\Q_p)^{\x}$ such that $u_pxu_p^{-1}=\si_p(x)$ for all $x\in F_p$.  Hence
$$
\CO_p=F_p\cap u_p^{-1}M(\Z_p)u_p.
$$

For $p\notin T$ set $u_p=1$ and let $\ut=\left(u_p\right)\in\A_{\fin}$.  By Lemma \ref{lem:Adeles}
$$
M(\A_{\fin})^{\x}=M(\hat{\Z})^{\x}\,M(\Q)^{\x}
$$
so there exists an element $u\in M(\Q)^{\x}$ such that for all primes $p$ we have
$
u^{-1}M(\Z_p)u=u_p^{-1}M(\Z_p)u_p.
$
This implies that
\begin{equation}
\label{eqn:padicid}
\CO_p=F_p\cap u^{-1}M(\Z_p)u
\end{equation}
for all primes $p$ and we deduce that
$
\CO=F\cap u^{-1}M(\Z)u=\CO_u.
$
This completes the proof of the lemma.
\qed

\subsection{Counting order embeddings}

Let $F/\Q$ be an algebraic number field and let $\CO$ be an order in $F$.  Let $I(\CO)$ \index{$I(\CO)$} be the set of all finitely generated $\CO$-submodules of $F$.  According to the Jordan-Zassenhaus Theorem \index{Jordan-Zassenhaus Theorem} (\cite{Reiner75}, Theorem 26.4), the set $[I(\CO)]$ of isomorphism classes of elements of $I(\CO)$ is finite.  Let $h(\CO)$ \index{$h(\CO)$} be the cardinality of the set $[I(\CO)]$, called the \emph{class number} \index{class number} of $\CO$.

Let $I$ be a non-trivial, finitely generated $\CO$-submodule of $F$.  Then, by \cite{Reiner75}, Theorem 10.6, $I$ is isomorphic as an $\CO$-module to an ideal of $\CO$ so also has the property that $I\ox\Q=F$.  If $J$ is another finitely generated $\CO$-submodule of $F$, such that $J\cong I$, then the isomorphism extends to an automorphism of $F$ considered as a module over itself, and hence $J=Ix$ for some $x\in F^{\x}$.  So $h(\CO)=\card([I(\CO)])=|I(\CO)/F^{\x}|$.

For a prime $p\in\Z$ let $F_p=F\ox\Q_p$\index{$F_p$}.  Let $\p_1,...\p_n$ be the prime ideals of $F$ above $p$.  By \cite{Neukirch99}, Chapter 1, Proposition 8.3 there exists a canonical $\Q_p$-algebra isomorphism defined by:
\begin{eqnarray}
\label{eqn:localproduct}
F_p=F\ox\Q_p & \cong   & \prod_{i=1}^n F_{\p_i} \\
    \al\ox x   & \mapsto & \left(\tau_i(\al).x\right)_{i=1}^n, \nonumber
\end{eqnarray}
where $F_{\p_i}$ denotes the completion of $F$ with respect to the $\p_i$-adic absolute value and $\tau_i$ the embedding of $F$ into $F_{\p_1}$.  If $\CO$ is an order in $F$ and $p\in\Z$ a prime then let $\CO_p=\CO\ox\Z_p$\index{$\CO_p$}.

Recall that $S$ is a finite, non-empty set of prime numbers with an even number of elements.  Let $C(S)$ \index{$C(S)$} be the set of all totally complex quartic fields $F$ such that all primes $p$ in $S$ are non-decomposed in $F$.  Note in particular that by the isomorphism (\ref{eqn:localproduct}) this means that for each field $F$ in $C(S)$ and for each $p\in S$ we have that $F_p=F\ox\Q_p$ is once again a field, namely the completion of $F$ at the unique place of $F$ above $p$.  For $F\in C(S)$ let $O_F(S)$ \index{$O_F(S)$} be the set of isomorphism classes of orders $\CO$ in $F$ which are maximal at all $p\in S$, ie. are such that the completion $\CO_p=\CO\ox\Z_p$ is the maximal order of the field $F_p$ for all $p\in S$.  Let $O(S)$ \index{$O(S)$} be the union of all $O_F(S)$, where $F$ ranges over $C(S)$.

Let $F\in C(S)$ and $\CO\in O_F(S)$, then, by Lemma \ref{lem:fieldembedding}, $F$ embeds into $M(\Q)$.  By Lemma \ref{lem:unit}, the group $M(\Z)^{\x}$ of units in $M(\Z)$ contains $\Ga$ as a subgroup.  By Lemma \ref{lem:orderembedding} we know that there is an embedding $\si$ of $F$ into $M(\Q)$ such that $\CO = \CO_{\si} = \si^{-1}\left(\si(F)\cap M(\Z)\right)$.  Let $u\in M(\Z)^{\x}$ and, as above, let $^u\si=u\si u^{-1}$ for embeddings $\si:F\ra M(\Q)$.  Then $\CO_{^u\si}=\CO_{\si}$, so the group $M(\Z)^{\x}$ acts on the set $\Si(\CO)$ \index{$\Si(\CO)$} of all embeddings $\si$ with $\CO_{\si}=\CO$.  This $M(\Z)^{\x}$ action also restricts to an action of $\Ga$ on $\Si(\CO)$.

We define \index{$\la_S(\CO)$}
$$
\la_S(\CO)=\la_S(F)=\prod_{p\in S} f_p(F),
$$
where $f_p(F)$ is the inertia degree of $p$ in $F$.  Our aim is to prove the following proposition.

\begin{proposition}
\label{pro:embeddings}
The quotient $M(\Z)^{\x}\bs\Si(\CO)$ is finite and has cardinality $h(\CO)\la_S(\CO)$.
\end{proposition}

We prepare for the proof of the propostition with the following discussion and then prove some lemmas, from which the proposition will follow.

Fix an embedding $\si:F\hra M(\Q)$ with $\CO=\CO_{\si}$ and consider $F$ as a subfield of $M(\Q)$ such that $\CO=F\cap M(\Z)$.  For $u\in M(\Q)^{\x}$ let
$
\CO_u =F\cap u^{-1}M(\Z)u.
$
Let $U$ be the set of all $u\in M(\Q)^{\x}$ such that
$
\CO=F\cap M(\Z)=F\cap u^{-1}M(\Z)u=\CO_u.
$
Let $a\in F^{\x}$, then
$$
F\cap a^{-1}u^{-1}M(\Z)ua=a^{-1}\left(F\cap u^{-1}M(\Z)u\right)a=F\cap u^{-1}M(\Z)u,
$$
so $F^{\x}$ acts on $U$ by multiplication from the right.  Let $v\in M(\Z)^{\x}$, then
$$
F\cap u^{-1}v^{-1}M(\Z)vu=F\cap u^{-1}M(\Z)u,
$$
so $M(\Z)^{\x}$ acts on $U$ by multiplication from the left.  Let $\tau\in\Si(\CO)$, so that $\CO_{\tau}=\CO$.  Then by the Skolem-Noether Theorem (Proposition \ref{pro:orders}(a)) there exists $u\in M(\Q)$ such that $\tau=\,^u\si=u\si u^{-1}$, ie.$\tau(F)=uFu^{-1}$.  Then 
\begin{eqnarray*}
\CO_u = F\cap u^{-1}M(\Z)u & = & u^{-1}\left(uFu^{-1}\cap M(\Z)\right)u \\
  & = & u^{-1}\left(\tau(F)\cap M(\Z)\right)u \\
  & = & \CO_{\tau} \= \CO.
\end{eqnarray*}
Conversely, it is clear that for all $u\in U$ we have $^u\si\in\Si(\CO)$, and that $^u\si=\ ^v\si$ if and only if $u=vx$ for some $x\in F^{\x}$.  Hence
$$
|M(\Z)^{\x}\bs U/F^{\x}|=|M(\Z)^{\x}\bs \Si(\CO)|,
$$
and the proposition will be proved if we can show that the left hand side equals $h(\CO)\la_S(\CO)$.

For $u\in U$ let
$$
I_u=F\cap M(\Z)u.
$$
Then $I_u$ is a finitely generated $\CO$-module in $F$.  We shall prove that the map \index{$\Psi$}
\begin{eqnarray*}
\Psi:M(\Z)^{\x}\bs U/F^{\x} & \ra & I(\CO)/F^{\x} \\
  u & \mapsto & I_u,
\end{eqnarray*}
is surjective and $\la_S(\CO)$ to one.  This will be done in the following lemmas through localisation and strong approximation.

Let $p\in\Z$ be a prime and let $I(\CO_p)$ be the set of all finitely generated $\CO_p$-submodules of $F_p$ and $[I(\CO_p)]$ the set of isomorphism classes.  Let $I_p$ be a non-trivial, finitely generated $\CO_p$-submodule of $F_p$.  Then, by \cite{Reiner75}, Theorem 10.6, the submodule $I_p$ is isomorphic as an $\CO_p$-module to an ideal of $\CO_p$ so also has the property that $I_p\ox\Q_p=F_p$.  If $J_p$ is another finitely generated $\CO_p$-submodule of $F_p$, such that $J_p\cong I_p$, then the isomorphism extends to an automorphism of $F_p$ considered as a module over itself, and hence $J_p=I_p x$ for some $x\in F_p^{\x}$.  So $\card([I(\CO_p)])=|I(\CO_p)/F_p^{\x}|$.

For a prime $p$ and $u_p\in M(\Q_p)^{\x}$ let 
$
\CO_{p,u_p}=F_p\cap u_p^{-1}M(\Z_p)u_p
$
and recall that
$
\CO_p=F_p\cap M(\Z_p).
$
Let $U_p$ be the set of all $u_p\in M(\Q_p)^{\x}$ such that
$
\CO_{p,u_p}=\CO_p.
$

As in the global case above, $F_p^{\x}$ acts on $U_p$ from the right and $M(\Z_p)^{\x}$ acts on $U_p$ from the left.  For $u_p\in U_p$ let
$$
I_{u_p}=F_p\cap M(\Z_p)u_p.
$$
Then $I_{u_p}$ is a finitely generated $\CO_p$-module in $F_p$.

\begin{lemma}
\label{lem:LocalIdeals}
Let $p$ be a prime.  Then $|I(\CO_p)/F_p^{\x}|=1$.
\end{lemma}
\prf
Suppose $p\in S$.  Then $p$ is non-decomposed in $F$ so there exists a unique prime ideal $\p$ of $F$ above $p$, so by the isomorphism (\ref{eqn:localproduct}) we can see that $F_p\cong F_{\p}$ and hence $F_p$ is itself a field.  By \cite{Neukirch99}, Chapter II, Proposition 3.9, every ideal of $\CO_p$ is principal and hence all non-zero ideals of $\CO_p$ are isomorphic.  The claim holds for $p\in S$.

Suppose $p\notin S$.  Let $\p_1,...,\p_n$ be the prime ideals of $F$ over $p$.  Under the isomorphism (\ref{eqn:localproduct}) an ideal in $\CO_p$ gives ideals in $\CO_{\p_1},...,\CO_{\p_n}$, where $\CO_{\p_1}$ denotes the ring of integers in the field $F_{\p_i}$.  As we saw above, these ideals are each unique up to isomorphism and hence the original ideal in $\CO_p$ is also unique up to isomorphism and the claim holds for $p\notin S$.
\qed

We now prove the surjectivity of $\Psi$.

\begin{lemma}
\label{lem:3}
The map $\Psi$ is surjective.
\end{lemma}
\prf
Let $I\subset\CO$ be an ideal and for a prime $p$ let $I_p=I\ox\Z_p$.  By Lemma \ref{lem:LocalIdeals}, for all primes $p$ the $\CO_p$-module $F_p\cap M(\Z_p)$ is isomorphic to $I_p$, so there exists $u_p\in F_p^{\x}$ such that
$
F_p\cap M(\Z_p)u_p=I_p.
$
We have further that
$
F_p\cap u_p^{-1}M(\Z_p)u_p  =  u_p^{-1}(F_p\cap M(\Z_p)u_p 
   =  F_p\cap M(\Z_p).
$

Let $\ut=(u_p)\in\A_{\fin}^{\x}$.  By Lemma \ref{lem:Adeles} there exist $\tilde{z}\in M(\hat{\Z})^{\x}$ and $u\in M(\Q)^{\x}$ such that $\ut=\tilde{z}u$.  Then for all primes $p$ we have that
$
F_p\cap u^{-1}M(\Z_p)u=F_p\cap M(\Z_p)
$
and
$
F_p\cap M(\Z_p)u=I_p.
$
We deduce from this that
$
F\cap u^{-1}M(\Z)u=F\cap M(\Z)
$
and
$
F\cap M(\Z)u=I.
$
This is what was required to complete the proof of the lemma.
\qed

In the following series of lemmas we prove that $\Psi$ is $\la_S(\CO)$ to one.

\begin{lemma}
\label{lem:Diag}
Let $K$ be a non-archimedian local field and $\CO_K$ its ring of integers.  For $n\in\N$ let $\D_n\subset\Mat_n(K)$ be the subspace of diagonal matrices and suppose that $u\in\GL_n(K)$ is such that
$$
\D_n\cap\Mat_n(\CO_K)\subset\D_n\cap u^{-1}\Mat_n(\CO_K)u.
$$
Then there exist $z\in\Mat_n(\CO_K)^{\x}$ and $x\in\D_n^{\x}$ such that $zu=x$.

If we suppose further that
$$
\D_n\cap\Mat_n(\CO_K)=\D_n\cap\Mat_n(\CO_K)u
$$
then we have $x\in(\D_n\cap\Mat_n(\CO_K))^{\x}$.
\end{lemma}
\prf
The first claim of the lemma will be proved by induction on $n$.  For $n=1$ we take $z=1$ and $x=u$ and the claim holds.  Assume now that the claim holds for some $n-1\in\N$.  For any ring $R$ let $\Upp_n(R)$ denote the subalgebra of upper triangular matrices in $\Mat_n(R)$.  We can choose an element $y\in\Mat_n(\CO_k)^{\x}$ such that $u'=yu\in\Upp_n(K)$ is upper triangular.  Then $u'$ has the form
$$
u'=\matrix{a}{b}{}{c},
$$
where $a\in K$, $b\in\Mat_{1\x(n-1)}(K)$ and $c\in\Upp_{n-1}(K)$.  Then
$$
(u')^{-1}=\matrix{a^{-1}}{-a^{-1}bc^{-1}}{}{c^{-1}}.
$$

Since $u$ satisfies
$
\D_n\cap\Mat_n(\CO_K)\subset\D_n\cap u^{-1}\Mat_n(\CO_K)u
$,
it follows that $u'$ satisfies
$
\D_n\cap\Upp_n(\CO_K)\subset\D_n\cap (u')^{-1}\Upp_n(\CO_K)u'.
$
It follows that there exist $\al\in\CO_K$, $\be\in\Mat_{1\x(n-1)}(\CO_K)$ and $\ga\in\Upp_{n-1}(\CO_K)$ such that
$$
\matrix{\al}{a^{-1}(\al b+\be c-bc^{-1}\ga c)}{}{c^{-1}\ga c}=(u')^{-1}\matrix{\al}{\be}{}{\ga}u'=\matrix{1}{}{}{0},
$$
where the matrix on the right is the $n$ by $n$ matrix with a one in the top left corner and all other entries zero.  It then follows that $\al=1$, $c^{-1}\ga c=0$ and $b=-\be c$.  Let
$$
w=\matrix{1}{\be}{}{I_{n-1}}\in\Mat_n(\CO_K)^{\x}.
$$
Then
$$
wu'=\matrix{1}{\be}{}{I_{n-1}}\matrix{a}{-\be c}{}{c}=\matrixtwo{a}{c}.
$$
so we have that
$$
wyu=\matrixtwo{a}{c}.
$$
It is then clear that $c$ satisfies
$$
\D_{n-1}\cap\Mat_{n-1}(\CO_K)\subset\D_{n-1}\cap c^{-1}\Mat_{n-1}(\CO_K)c
$$
so by the inductive hypothesis there exist $z'\in\Mat_{n-1}(\CO_K)^{\x}$ and $x'\in\D_{n-1}^{\x}$ such that $z'c=x'$.  Setting
$$
z=\matrixtwo{1}{z'}wy\in\Mat_n(\CO_K)^{\x}\ \ \ \textrm{and}\ \ \ x=\matrixtwo{a}{x'}\in\D_n^{\x}
$$
we get that $zu=x$ and the proof of the first claim is complete.

If we also assume that $u$ satisfies
$$
\D_n\cap\Mat_n(\CO_K)=\D_n\cap\Mat_n(\CO_K)u
$$
then since $z\in\Mat_n(\CO_K)^{\x}$ and $zu=x\in\D_n^{\x}$ it follows that
$$
\D_n\cap\Mat_n(\CO_K)=\D_n\cap\Mat_n(\CO_K)zu=(\D_n\cap\Mat_n(\CO_K))x.
$$
This can be the case only if $x\in(\D_n\cap\Mat_n(\CO_K))^{\x}$ so the lemma is proved.
\qed

\begin{lemma}
\label{lem:1}
For $p\notin S$ we have $|M(\Z_p)^{\x}\bs U_p/F_p^{\x}|=1$.
\end{lemma}
\prf
Let $p\notin S$ and let $u_p\in U_p$ so that
\begin{equation}
\label{eqn:Fp1}
F_p\cap M(\Z_p)=F_p\cap u_p^{-1}M(\Z_p)u_p.
\end{equation}
By Lemma \ref{lem:LocalIdeals}, there exists $\la_p\in F_p^{\x}$ such that
$
F_p\cap M(\Z_p)=F_p\cap M(\Z_p)u_p\la_p.
$
Replacing $u_p$ with $u_p\la_p$ we may assume that $u_p$ satisfies
\begin{equation}
\label{eqn:Fp2}
F_p\cap M(\Z_p)=F_p\cap M(\Z_p)u_p.
\end{equation}
We shall show that there exist $z_p\in M(\Z_p)^{\x}$ and $x_p\in F_p^{\x}$ such that $z_pu_px_p$ is the identity element in $M(\Q_p)$.

Let $\ep$ be an integral element of $F$ generating $F$ over $\Q$, and let $L$ be an extension of $\Q$ containing all the zeros of the minimal polynomial of $\ep$.  Let $L_p=L\ox\Q_p$, note that neither this nor $F_p$ are necessarily fields, rather a finite direct product of fields.  The embedding $F\hra M(\Q)$ gives us an embedding $F_p\hra M(\Q_p)$.  Since $p\notin S$ we have $M(\Q_p)\cong\Mat_4(\Q_p)$ so we get an embedding
$$
F_p\ox L\hra M(\Q_p)\ox L\cong\Mat_4(L_p).
$$

Let $\CO_L$ be the ring of integers in $L$ and let $\CO_{L_p}=\CO_L\ox\Z_p$.  Let $\D$ denote the subspace of diagonal matrices in $\Mat_4(L_p)$.  It follows from our choice of $L$ that $\ep$ is diagonalisable in $M(L_p)$ by an element $A\in M(\CO_{L_p})^{\x}$.  Since $F_p=\Q_p(\ep)$, the whole of $F_p\ox L$ is simultaneously diagonalisable in $M(L_p)$, that is, we have that $A^{-1}(F_p\ox L)A$ is a subspace of $\D$, and by dimensional considerations we must have that $A^{-1}(F_p\ox L)A=\D$.  Then, tensoring (\ref{eqn:Fp1}) and (\ref{eqn:Fp2}) with $\CO_L$ and conjugating by $A$ we get
\begin{equation}
\label{eqn:Lp1}
\D\cap M(\CO_{L_p})=\D\cap\bar{u}_p^{-1}M(\CO_{L_p})\bar{u}_p
\end{equation}
and
\begin{equation}
\label{eqn:Lp2}
\D\cap M(\CO_{L_p})=\D\cap M(\CO_{L_p})\bar{u}_p,
\end{equation}
where $\bar{u}_p=u_p A$.

We are currently working in the matrix algebra $M(L_p)\cong\Mat_4(L_p)$.  Let $\p_1,...,\p_n$ be the prime ideals of $L$ above $p$.  By virtue of the isomorphism (\ref{eqn:localproduct}) we can consider separately each $\p_i$-adic component of the entries of the matrices.  It then follows from Lemma \ref{lem:Diag} and equations (\ref{eqn:Lp1}) and (\ref{eqn:Lp2}) that there exist $\bar{z}\in M(\CO_{L_p})^{\x}$ and $\bar{x}\in(\D\cap M(\CO_{L_p}))^{\x}$ such that $\bar{z}\bar{u}_p=\bar{x}$.  Setting
$
z=A\bar{z}\in M(\CO_{L_p})^{\x}
$
and
$
x=A\bar{x}A^{-1}\in ((F_p\ox L)\cap M(\CO_{L_p}))^{\x}
$
we get that $zu_p=x$.

Consider now the trace map $\tr_{L/\Q}$ of the field extension $L/\Q$.  The image of $\CO_L$ under this map is an ideal in $\Z$.  Let $\nu\in\Z$ be the greatest such that $\tr_{L/\Q}(\CO_L)\subset p^{\nu}\Z$.  We define a linear map
\begin{eqnarray*}
T:M(L_p)=M(\Q_p)\ox L & \ra &     M(\Q_p)\ox\Q\cong M(\Q_p) \\
           x\ox y     & \mapsto & x\ox p^{-\nu}(\tr_{L/\Q}\,y).
\end{eqnarray*}
Note that $T(M(\CO_{L_p}))=T(M(\Z_p)\ox\CO_L)\subset M(\Z_p)\ox\Z\cong M(\Z_p)$.

We denote by $\CO_{L_1}$ the ring of integers of the field $L_{\p_1}$.  Let $\xi\in(F\ox L_{\p_1})\cap M(\CO_{L_1})$ and let $x_1$ and $z_1$ be the $\p_1$ components of $x$ and $z$.  Then $z_1u_p=x_1$.  Now set $\al=x_1^{-1}\xi\in (F\ox L_{\p_1})\cap M(\CO_{L_1})$.  From (\ref{eqn:Fp1}) we deduce that
$$
(F\ox L_{\p_1})\cap M(\CO_{L_1})=(F\ox L_{\p_1})\cap u_p^{-1}M(\CO_{L_1})u_p
$$
so there exists $\be\in M(\CO_{L_1})$ such that $\al=u_p^{-1}\be u_p$, or equivalently $u_p\al=\be u_p$.  Setting $\ze=z_1\be\in M(\CO_{L_1})$, we then get
$$
\ze u_p=z_1\be u_p=z_1u_p\al=x_1\al=\xi.
$$
Hence we can see that $\xi u_p^{-1}\in M(\CO_{L_1})$.

We now consider the linear map defined by
\begin{eqnarray*}
T':M(\Q)\ox L_{\p_1} & \ra &     M(\Q)\ox\Q_p=M(\Q_p) \\
      x \ox y        & \mapsto & x\ox p^{-\nu}(\tr_{L_{\p_1}/\Q_p}\,y),
\end{eqnarray*}
where $\tr_{L_{\p_1}/\Q_p}$ is the trace map of the field extension $L_{\p_1}/\Q_p$.  We note that $T'$ maps $(F\ox L_{\p_1})\cap M(\CO_{L_{\p}})$ to $F_p\cap M(\Z_p)$.  The map $T'$ induces a  linear map
$$
T'':\frac{M(\CO_{L_1})}{\p_1M(\CO_{L_1})}\ra\frac{M(\Z_p)}{pM(\Z_p)}.
$$
Note that $M(\CO_{L_1})/\p_1M(\CO_{L_1})$ is a vector space over the field $\CO_{L_1}/\p_1\CO_{L_1}$ and $M(\Z_p)/pM(\Z_p)$ is a vector space over the field $\Z_p/p\Z_p\cong\F_p$.  By the choice of the integer $\nu$ we can see that the map $T''$ is not the zero map, that is, its kernel is a proper subspace of $M(\CO_{L_1})/\p_1M(\CO_{L_1})$.  An element of $M(\CO_{L_1})/\p_1M(\CO_{L_1})$ not in the kernel of $T''$ corresponds to an element of $M(\CO_{L_1})$ whose image under $T'$ is in $M(\Z_p)^{\x}$.  We choose an element $\xi\in(F\ox L_{\p_1})\cap M(\CO_{L_1})$ such that $T'(\xi)\in (F_p\cap M(\Z_p))^{\x}$ and $T'(\xi u_p^{-1})\in M(\Z_p)^{\x}$.

We shall write $\xi$ also for the element of $(F_p\ox L)\cap M(\CO_{L_p})$ with $\p_1$ component equal to $\xi$ and all other components zero.  Similarly we denote by $\ze$ the element of $M(\CO_{L_p})$ with $\p_1$ component equal to $\xi u_p^{-1}$ and all other components zero.  Then $\ze u_p=\xi$.  From the definitions of the maps $T$ and $T'$ and the isomorphism (\ref{eqn:localproduct}) it follows that
$
T(\xi)=T'(\xi')\in(F_p\cap M(\Z_p))^{\x}
$
and
$
T(\ze)=T'(\ze)\in M(\Z_p)^{\x}.
$

Now set $x_p=T(\xi)\in(F_p\cap M(\Z_p))^{\x}$ and $z_p=T(\ze)\in M(\Z_p)^{\x}$.  We then have that
$
z_pu_p=x_p
$
and the lemma is proved.
\qed

\begin{lemma}
\label{lem:2}
For $p\in S$ we have $|M(\Z_p)^{\x}\bs U_p/F_p^{\x}|=f_p(F)$.
\end{lemma}
\prf
Let $p\in S$.  We then have that $M(\Q_p)$ is a division algebra so in particular $M(\Q_p)^{\x}=M(\Q_p)\smallsetminus\{ 0\}$.  Further, since $\Q_p$ is a local field, $M(\Z_p)$ is the unique maximal order of $M(\Q_p)$ (cf.Proposition \ref{pro:orders}(d)) and hence $u_p^{-1}M(\Z_p)u_p=M(\Z_p)$ for all $u_p\in M(\Q_p)^{\x}$.  Thus the claim is equivalent to
$$
|M(\Z_p)^{\x}\bs M(\Q_p)^{\x}/F_p^{\x}|=f_p(F).
$$
Let $v_M$ be the $p$-adic valuation on $M(\Q_p)$ and let
$$
e_M=e(M(\Q_p)/\Q_p)=|v_M(M(\Q_p)^{\x})/v_M(\Q_p^{\x})|
$$
be the ramification index of $\Q_p$ in $M(\Q_p)$.  The valuation $v_M$ is a surjective group homomorphism of $M(\Q_p)^{\x}$ onto $\Z$ with kernel $M(\Z_p)^{\x}$.  It follows that
$$
e_M=|M(\Z_p)^{\x}\bs M(\Q_p)^{\x}/\Q_p^{\x}|.
$$
By Proposition 17.7 of \cite{Pierce82} we have
$$
|M(\Z_p)^{\x}\bs M(\Q_p)^{\x}/\Q_p^{\x}|=e_M=\deg M(\Q_p)=4.
$$

Let $e_p(F)$ be the ramification index of $p$ in $F$.  Then by \cite{Neukirch99}, Chapter I, Proposition 8.2, and since $p$ does not split in $F$, we have
\begin{equation}
\label{eqn:fundid}
e_p(F)f_p(F)=[F:\Q]=4.
\end{equation}
Let $\CCO$ be the valuation ring of $F_p$ and let $\Pi=\pi\CCO$ be the unique maximal ideal of $\CCO$, where $\pi\in\CCO$ is a generator of the principal ideal $\Pi$.  Let $\p=p\Z_p$ be the unique maximal ideal of $\Z_p$.  Then $\p=\Pi\cap\Z_p$ and $p=\pi^{e_p(F)}$.  Thus $v_M(\pi)=v_M(p)/e_p(F)$, which implies that
$$
|M(\Z_p)^{\x}\bs M(\Q_p)^{\x}/F_p^{\x}|=\frac{|M(\Z_p)^{\x}\bs M(\Q_p)^{\x}/\Q_p^{\x}|}{e_p(F)}=\frac{4}{e_p(F)},
$$
and hence from (\ref{eqn:fundid}) we get
$$
|M(\Z_p)^{\x}\bs M(\Q_p)^{\x}/F_p^{\x}|=f_p(F),
$$
as required.
\qed

\begin{lemma}
$\Psi$ is $\la_S(\CO)$ to one.
\end{lemma}
\prf
Recall that $\CO\in O_F(S)$ for some $F\in C(S)$.  Let $\hat{F}=\prod'_p F_p$ be the restricted product over all primes $p$, that is, all but finitely many components are integral, and let $\hat{U}$ be the set of all $\ut=(u_p)\in M(\A_{\fin})^{\x}$ such that $\CO_p=F_p\cap u_p^{-1}M(\Z_p)u_p=\CO_{p,u_p}$ for all $p$.

For $I\in I(\CO)$ let $U_I=\Psi^{-1}(IF^{\x})$.  Let $u\in U_I$ for some $I\in I(\CO)$.  We denote by $\hat{u}$ the element of $M(\A_{\fin})$ each of whose components is $u$.  Setting $u_p=u$ for all primes $p$, we see that for all $p$
$$
\CO_p=\CO\ox\Z_p=\CO_u\ox\Z_p=\CO_{p,u_p}
$$
so $\hat{u}\in\hat{U}$.

For each $I\in I(\CO)$ we define the map
$$
\Phi_I:M(\Z)^{\x}\bs U_I/F^{\x}\ra M(\hat{\Z})^{\x}\bs\hat{U}/\hat{F}^{\x}
$$
by setting $\Phi_I(M(\Z)^{\x}uF^{\x})=M(\hat{\Z})^{\x}\hat{u}\hat{F}^{\x}$.  We shall show that $\Phi_I$ is a bijection.  Since from Lemmas \ref{lem:1} and \ref{lem:2} we know that $|M(\hat{\Z})^{\x}\bs\hat{U}/\hat{F}^{\x}|=\la_S(\CO)$, the lemma will then follow.

First we show surjectivity.  Let $\ut=(u_p)\in\hat{U}$, so that $\CO_{p,u_p}=\CO_p$ for all primes $p$.  For all primes $p$, by Lemma \ref{lem:LocalIdeals} we have $I\ox\Q_p\cong F_p\cap M(\Z_p)u_p$ so there exists $x_p\in F_p^{\x}$ such that $I\ox\Q_p=F_p\cap M(\Z_p)u_px_p$.  We replace $u_p$ with $u_px_p$.  By Lemma \ref{lem:Adeles} there exist $\tilde{z}=(z_p)\in M(\hat{\Z})^{\x}$ and $u\in M(\Q)^{\x}$ such that $\ut=\tilde{z}\hat{u}$.  It follows that $\CO_{p,u}=\CO_p$ for all primes $p$ and hence we get that $\CO_u=\CO$, so that $u\in U$.  Furthermore, since
$$
I\ox\Q_p=F_p\cap M(\Z_p)u_p=F_p\cap M(\Z_p)u
$$
for all primes $p$, it follows that $I=F\cap M(\Z)u$, so $u\in U_I$.  We have proven that $\Phi_I$ is surjective.

For the proof of injectivity, suppose that $u,v\in U_I$ are such that $\Phi_I(u)=\Phi_I(v)$.  The assumption that $u,v\in U_I$ means that there is an isomorphism of $\CO$-modules $F\cap M(\Z)u\cong F\cap M(\Z)v$.  As we saw above this implies that there exists $x\in F^{\x}$ such that $F\cap M(\Z)u=F\cap M(\Z)vx$.  Replacing $v$ with $vx$ it follows that
\begin{equation}
\label{eqn:IdealEq}
\hat{F}\cap M(\hat{\Z})\hat{u}=\hat{F}\cap M(\hat{\Z})\hat{v}.
\end{equation}
The assumption that $\Phi_I(u)=\Phi_I(v)$ means that there exist $\tilde{z}\in M(\hat{\Z})^{\x}$ and $\tilde{x}\in\hat{F}^{\x}$ such that $\hat{u}=\tilde{z}\hat{v}\tilde{x}$.  Together with (\ref{eqn:IdealEq}) this implies that
$$
\hat{F}\cap M(\hat{\Z})\hat{v}=\hat{F}\cap M(\hat{\Z})\hat{v}\tilde{x}
$$
and hence $\tilde{x}\in\hat{v}^{-1}M(\hat{\Z})^{\x}\hat{v}$.  In other words, there exists $\tilde{y}\in M(\hat{\Z})^{\x}$ such that $\hat{v}\tilde{x}=\tilde{y}\hat{v}$.  We deduce that
$$
\hat{u}=\tilde{z}\hat{v}\tilde{x}=(\tilde{z}\tilde{y})\hat{v}.
$$
Replacing $\tilde{z}$ with $\tilde{z}\tilde{y}\in M(\hat{\Z})^{\x}$, we now have that
$$
\tilde{z}=\hat{u}\hat{v}^{-1}\in M(\hat{\Z})^{\x}\cap M(\Q)^{\x}=M(\Z)^{\x}.
$$
Writing $z=\tilde{z}\in M(\Z)^{\x}$ we get that $u=zv$ and injectivity is proven.
\qed
In this final section we complete the proof of the Main Theorem.  We do this by using our results on orders in subfields of division algebras of degree four to interpret the Prime Geodesic Theorem algebraically.  This will then yield information about the class numbers of those orders.

\section{Comparing Geodesics and Orders}
    \label{ch:Compare}
\subsection{Regular geodesics}
Let $F$ be a number field with $r$ real embeddings and $s$ pairs of complex conjugate embeddings and let $t=r+s-1$.  Let $\CO_F$ denote the ring of integers of $F$ and $\CO_F^{\x}$ the unit group therein.  By Dirichlet's unit theorem (\cite{Neukirch99}, I Theorem 7.4), there exist units $\ep_1,...,\ep_t\in\CO_F^{\x}$ such that every $\ep\in\CO_F^{\x}$ can be written uniquely as a product
$$
\ep=\zeta\ep_1^{\nu_1}\cdots\ep_t^{\nu_t},
$$
where $\zeta$ is a root of unity and $\nu_i\in\Z$.  The set $\{\ep_1,...,\ep_t\}$ is not uniquely determined.  We shall refer to such a set of units as a \emph{system of fundamental units}\index{fundamental unit}.  If $\CO\subset\CO_F$ is an order of $F$ and $\{\ep_1,...,\ep_t\}$ a system of fundamental units then there exist $\mu_i\in\N$ and units $\bar{\ep}_i=\ep_i^{\mu_i}$ such that every $\ep\in\CO^{\x}$ can be written uniquely as a product
$$
\ep=\zeta\bar{\ep}_1^{\nu_1}\cdots\bar{\ep}_t^{\nu_t}
$$
where $\zeta$ is a root of unity and $\nu_i\in\Z$.  Once again, the set $\{\bar{\ep}_1,...,\bar{\ep}_t\}$ is not uniquely determined.  We shall call such a set of units \emph{a system of fundamental units of $\CO$}.

If $\tau_1,...,\tau_{t+1}$ are distinct embeddings of $F$ into $\C$ which are pairwise non-conjugate, then the \emph{regulator} \index{regulator} $R(\CO)$ of an order $\CO\subset\CO_F$ is defined as the absolute value of the determinant of an arbitrary minor of rank $t$ of the matrix
$$
\left(\log|\tau_i(\bar{\ep}_j)|\right)_{ij},
$$
where $1\leq i\leq t+1$ and $1\leq j\leq t$.  The value of the regulator does not depend on the choice of a system of fundamental units.

We shall sometimes refer, somewhat imprecisely, to the elements of a system of fundamental units simply as \emph{fundamental units}.  In the case $t=1$, let $\ep_1$ be a fundamental unit of an order $\CO$ in $F$.  We shall refer to the elements of the set $\{\zeta\ep_1^{\pm 1}:\zeta\textrm{ root of unity in }F\}$ as \emph{the fundamental units of $\CO$}, or in the case that $\CO=\CO_F$ \emph{the fundamental units of $F$}.

In the case that $F$ is a totally complex quartic field we have $t=1$, so we can choose a fundamental unit $\ep_1$ such that for all $\ep\in\CO_F^{\x}$ we have $\ep=\zeta\ep_1^{\nu_1}$, where $\zeta$ is a root of unity and $\nu_1\in\Z$.  Then $\ep_1$ is determined up to inversion and multiplication by a root of unity.  In this case the regulator of an order $\CO\subset\CO_F$ of $F$ with fundamental unit $\bar{\ep}_1$ satisfies
$$
R(\CO)=\left|\,\log|\tau(\bar{\ep}_1)|\,\right|,
$$
where $\tau$ is any embedding of $F$ into $\C$.

We consider subfields of the field $F$.  Apart from $\Q$ these are all quadratic.  By Dirichlet's unit theorem, real quadratic fields contain a single fundamental unit (up to inversion and multiplication by $\pm 1$) and complex quadratic fields contain no fundamental unit.  We can see then that $F$ can have at most one real subfield.  Let $C^r(S)$ \index{$C^r(S)$} and $C^c(S)$ \index{$C^c(S)$} be respectively the subsets of $C(S)$ consisting of fields with and without a real quadratic subfield.  Let $O^r(S)$ \index{$O^r(S)$} and $O^c(S)$ \index{$O^c(S)$} be respectively the subsets of $O(S)$ consisting of orders in fields in $C^r(S)$ and $C^c(S)$.  Then $C(S)=C^r(S)\cup C^c(S)$ and $O(S)=O^r(S)\cup O^c(S)$.

Let $[\ga]\in\CE_P^{p,\reg}(\Ga)$.  Then by Lemma \ref{lem:centralisers} the centraliser $M(\Q)_{\ga}$ is a totally complex quartic field which embeds into $M(\Q)$ and we denote this subfield of $M(\Q)$ by $F_{\ga}$.  We say that an element $\ga\in\Ga$ is weakly neat if the adjoint $\Ad(\ga)$ has no non-trivial roots of unity as eigenvalues.  Let $\CE_P^{p,\rw}(\Ga)$ \index{$\CE_P^{p,\rw}(\Ga)$} be the set of primitive, regular, weakly neat elements in $\CE_P(\Ga)$.

\begin{lemma}
\label{lem:WeaklyNeat}
Let $[\ga]\in\CE_P^{p,\reg}(\Ga)$.  Then $F_{\ga}\in C^c(S)$ if and only if $[\ga]\in\CE_P^{p,\rw}(\Ga)$.
\end{lemma}
\prf
Let $\CO$ be the maximal order in $F$ and let $\CO^{\x}$ denote the group of units.  By Dirichlet's unit theorem we can see that $F_{\ga}$ contains a real quadratic subfield if and only if $\CO^{\x}\cap\R\neq\{\pm 1\}$.  Let $\ga_1\in\CO^{\x}$ be a fundamental unit.  The element $\ga$ is a unit $F_{\ga}$ so there exists a root of unity $\ze$ and $n\in\Z$ such that $\ga=\ze\ga_1^n$.  If there also exist a root of unity $\xi$ and $m\in\Z\smallsetminus\{0\}$ such that $\xi\ga_1^m\in\R$ then $\ga_1$ is not weakly neat and so neither is $\ga$.

Conversely, suppose $\ga$ is not weakly neat.  We know that $\ga$ is conjugate in $G$ to an element $a_{\ga}b_{\ga}\in A^-B$.  From the assumption that $\ga$ is not weakly neat it follows that there exists $m\in\N$ such that one of the components of $b_{\ga}^m$ is diagonal.  Since the group $AB$ is abelian we have that $\ga^m$ is conjugate in $G$ to $a_{\ga}^m b_{\ga}^m$.  Since $\ga^m\in\CE_P(\Ga)$ it follows from Lemma \ref{lem:centralisers} that $b_{\ga}^m$ is diagonal and hence $\ga^m$ generates a real quadratic subfield of $F_{\ga}$.  The lemma follows.
\qed

Define the map
$$
\th:\CE_P^{p,\rw}(\Ga)\ra O^c(S)
$$
by
$$
\th:[\ga]\mapsto F_{\ga}\cap M(\Z).
$$
Let $\ga,\d\in\Ga$.  Then $F_{\d\ga\d^{-1}}=\d F_{\ga}\d^{-1}$, and since by Lemma \ref{lem:unit} the group $\Ga$ is a subgroup of the multiplicative group $M(\Z)^{\x}$, we can conclude that
$$
F_{\d\ga\d^{-1}}\cap M(\Z)=\d\left(F_{\ga}\cap M(\Z)\right)\d^{-1}\cong F_{\ga}\cap M(\Z).
$$
Hence the map is well defined.

\begin{proposition}
The map $\th$ is surjective.
\end{proposition}

\prf
Let $\CO\in O^c(S)$ be an order in the field $F\in C^c(S)$.  Since $\deg F=4$, Proposition \ref{pro:subfields}(b) implies that $F$ is a maximal subfield of $M(\Q)$.

Lemma \ref{lem:orderembedding} tells us that there exists an embedding $\si:F\ra M(\Q)$ such that
$$
\CO=\CO_{\si}=\si^{-1}\left(\si(F)\cap M(\Z)\right).
$$
Let $\ep$ be a fundamental unit in $\CO^{\x}$.  Let $\ga$ be the image of $\ep$ under the embedding $\si$.  The unit $\ep$ is integral in $F$, so $\ga\in M(\Z)$.

First we show that $\ga\sim_G a_{\ga}b_{\ga}$ for some $a_{\ga}\in A^{-}$, $b_{\ga}\in B$.  The field $F$ has two pairs of conjugate embeddings into $\C$.  Let $\si_1$, $\si_2$ be two distinct, non-conjugate embeddings of $F$ into $\C$.  The $\R$-algebra $F\ox_{\Q}\R$ is isomorphic to the commutative $\R$-algebra $\C\oplus\C$ via the map 
$$
\al\ox x \mapsto \left(x\si_1(\al),x\si_2(\al)\right).
$$
We note that $\CO\ox_{\Z}\R =F\ox_{\Q}\R$ and we get the series of inclusions
$$
\CO \subset \CO\ox_{\Z}\R = F\ox_{\Q}\R \subset M(\Q)\ox_{\Q}\R = \Mat_4(\R) \subset \Mat_4(\C).
$$
Thus we see that $\ga$ may be considered as an element, which we will denote $\bar{\ga}$, of a commutative subalgebra of $\Mat_4(\C)$ isomorphic to $\C\oplus\C$.  The matrix $\bar{\ga}$ has real entries and can be diagonalised in $\Mat_4(\C)$ to the matrix
$$
X=\matrixfour{\si_1(\ga)}{\overline{\si_1(\ga)}}{\si_2(\ga)}{\overline{\si_2(\ga)}},
$$
since the eigenvalues of $\bar{\ga}$ are the roots of its minimal polynomial, ie. the values of $\ga$ under its different embeddings into $\C$.  Let $a,b\in\R^+$; $\th,\phi\in [-\pi,\pi]$ be such that $\si_1(\ga)=ae^{i\th}$ and $\si_2(\ga)=be^{i\phi}$.  Since $\ga\in M(\Z)$, we have $\prod_{\si}\si(\ga)=N_{F|\Q}(\ga)=\det\ga=\pm1$, where $N_{F|\Q}$ is the field norm and the product is over all embeddings $\si$ of $F$ into $\C$ (\cite{Neukirch99}, I Proposition 2.6(iii)).  Hence $a^2 b^2=\pm 1$.  But $a$ and $b$ are both real so we must have $a^2 b^2=1$, and so $b=a^{-1}$ and we have
$$
X=\matrixfour{ae^{i\th}}{ae^{-i\th}}{\rez{a}e^{i\phi}}{\rez{a}e^{-i\phi}}.
$$
Without loss of generality we assume $a<b$ so that $a\in(0,1)$.
Let $R$ be the $2\x 2$ matrix
$$
\rez{\sqrt{2}}\matrix{i}{-1}{1}{-i}
$$
and $Y$ the $4\x 4$ matrix
$$
\matrixtwo{R}{R}\in\Mat_4(\C)
$$
Then $X$ is conjugate in $\Mat_4(\C)$ via $Y$ to the real matrix
$$
X'=\matrixtwo{aR(\th)}{\rez{a}R(\phi)}\in A^-B.
$$
So $\bar{\ga}$ is conjugate in $\Mat_4(\C)$ to an element of $A^-B$, hence by \cite{Lang02}, XIV Corollary 2.3, the matrix $\bar{\ga}$ is conjugate in $\Mat_4(\R)$ to $X'$, that is $(Z')^{-1}\bar{\ga}Z'=X'$ for some $Z'\in\Mat_4(\R)$.  Since $Z'$ is invertible we have that $n=\det\,Z'\neq 0$.  Let $Z=|n|^{-\frac{1}{4}}Z'$.  Then $\det\,Z=\pm 1$ and $Z^{-1}\bar{\ga}Z=X'$, so that $\ga$ is conjugate in $\SL_4^{\pm}(\R)$ to an element of $A^-B$.  Suppose that $\det\,Z=-1$.  Let $W$ be the matrix
$$
W=\left( \begin{array}{cccc} 
	       0 & 1 & \ & \ \\ 
	       1 & 0 & \ & \ \\  
	       \ & \ & 1 & 0 \\  
	       \ & \ & 0 & 1 \end{array}\right).
$$
Then $ZW\in G$ and $(ZW)^{-1}\bar{\ga}ZW=X''$, where
$$
X''=\matrixtwo{aR(-\th)}{\rez{a}R(\phi)}\in A^-B.
$$
Hence $\ga$ is conjugate in $G$ to an element of $A^-B$.

We also saw above that $\det\ga=1$ and so $\ga\in\Ga$ and hence $\ga\in\CE_P(\Ga)$.  It then follows from Lemma \ref{lem:centralisers} that $\ga$, and hence also $\ep$, generates either a real quadratic or a totally complex quartic field over $\Q$.  The field generated by $\ep$ over $\Q$ is a subfield of $F$.  However $F\in C^c(S)$ so has no real quadratic subfields.  Hence $\Q(\ga)\cong\Q(\ep)$ is a totally complex quartic field and Lemma \ref{lem:centralisers} tells us that $\ga$ is regular.  It then follows from Lemma \ref{lem:WeaklyNeat} that $\ga$ is weakly neat.

It remains to show that $\ga$ is primitive in $\Ga$.  Note that $\Q(\ep)$ is a quartic field contained in and hence equal to the field $F$.  Recall that we write $F_{\ga}=M(\Q)_{\ga}$ for the centraliser of $\ga$ in $M(\Q)$.  Then we have that $F_{\ga}=\Q(\ga)\cong F$ and $\CO^{\x}\cong F_{\ga}\cap M(\Z)^{\x}$.  The fundamental unit $\ep$ is primitive in $\CO^{\x}$ and so $\ga=\si(\ep)$ is primitive in $F_{\ga}\cap M(\Z)^{\x}$.  Suppose $\ga$ is not primitive in $\Ga$.  Then there exists $\ga_0\in\Ga\subset M(\Z)^{\x}$ such that $\ga=\ga_0^{\kappa}$ for some $\kappa\in\N$.  However we also have $\ga_0\in F_{\ga}$ and so $\ga_0\in F_{\ga}\cap M(\Z)^{\x}$, which is a contradiction.  Hence $\ga$ is primitive in $\Ga$ and we have $\ga\in\CE_P^{p,\reg}(\Ga)$.

This concludes the proof of the surjectivity of $\th$.
\qed

\begin{lemma}
\label{lem:rootsOfUnity}
Let $\ze$ be a torsion element of $M(\Z)$.  Then $\det(\ze)=1$.
\end{lemma}
\prf
By Proposition \ref{pro:subfields}(b), the element $\ze$ generates a subfield $\Q(\ze)$ of $M(\Q)$ of degree 1, 2 or 4 over $\Q$.  If $[\Q(\ze):\Q]=1$ or 2 then by Lemma \ref{lem:unit} we have $\det(\ze)=1$.

For a primitive $n^{\rm th}$ root of unity $\xi$, the degree $[\Q(\xi):\Q]$ is equal to $\varphi(n)=\#\{1\leq m\leq n:(m,n)=1\}$.  From \cite{Chandrasekharan68}, Chapter II, \S 2 we see that
$$
\varphi(n)=\prod_{i=1}^r p_i^{a_i-1}(p_i-1),
$$
where $n=p_1^{a_1}...p_r^{a_r}$ is the prime decomposition of $n$.  It is then straightfoward to show that $[\Q(\xi):\Q]=4$ if and only if $n=5,8,10$ or 12.  In these cases the minimal polynomial of $\xi$ is respectively $x^4+x^3+x^2+x+1$, $x^4+1$, $x^4-x^3+x^2-x+1$ or $x^4-x^2+1$.  In each of these cases we read off from the constant term of the minimal polynomial that $N_{\Q(\xi)|\Q}(\xi)=1$.  The lemma follows.
\qed

For $\CO\in O(S)$ let $\mu_{\CO}$ be the number of roots of unity in the order $\CO$.  For $[\ga]\in\CE_P^{p,\reg}(\Ga)$ let $\mu_{\ga}$ denote the number of roots of unity in the order $F_{\ga}\cap M(\Z)$.  Note that $F_{\ga}\cap M(\Z)$ is equal to the centraliser $M(\Z)_{\ga}$ of $\ga$ in $M(\Z)$, and by Lemma \ref{lem:rootsOfUnity} the roots of unity in $M(\Z)_{\ga}$ are all in $\Ga$ and hence in $\Ga_{\ga}$.  It follows that $\mu_{\ga}$ is equal to the cardinality of the torsion part of $\Ga_{\ga}$.  It is also immediate from the definitions that, writing $\CO_{\ga}=F_{\ga}\cap M(\Z)$, we have $\mu_{\ga}=\mu_{\CO_{\ga}}$.

\begin{lemma}
\label{lem:SplitFactors}
Let $\CO\in O(S)$ be an order in the field $F\in C(S)$ and let $f(x)$ be the minimal polynomial over $\Q$ of a fundamental unit in $\CO$.  Then the number of roots of $f(x)$ in $F$ which are also fundamental units in $\CO$ is independent of the choice of fundamental unit.  We call this number $\ka(\CO)$\index{$\ka(\CO)$}.  Note that $\ka(\CO)$ is equal to 1,2 or 4.
\end{lemma}

\prf
Let $\ep$ be a fundamental unit in $\CO$ with minimal polynomial $f(x)$ over $\Q$.  Suppose there exists another fundamental unit $\d\in\CO^{\x}$ with $\d\neq\ep$ which is also a root of $f(x)$.  Then there exist embeddings $\al$ and $\be$ of $F$ into $\C$ such that $\al(\ep)=\be(\d)$.

Let $\eta\in\CO^{\x}$ also be a fundamental unit in $\CO$ with minimal polynomial $g(x)$ over $\Q$.  Then $\eta=\ze\ep^{\pm 1}$ for some root of unity $\ze\in\CO^{\x}$.  We have
$$
\al(\eta)=\al(\ze\ep^{\pm 1})=\al(\ze)\al(\ep)^{\pm 1}=\al(\ze)\be(\d)^{\pm 1}.
$$
Let $n$ be the order of the root of unity $\ze$.  The order $\CO$ is a ring so contains all roots of unity of order $n$ and hence there exists an $n^{th}$ root of unity $\xi\in\CO^{\x}$ such that $\be(\xi)=\al(\ze)$.  Let $\th$ be the fundamental unit $\xi\d^{\pm 1}\in\CO^{\x}$.  Then $\al(\eta)=\be(\th)$ and so $\th$ is also a root of $g(x)$.  The lemma follows.
\qed

\begin{proposition}
\label{pro:ThetahlaToOne}
The map $\th$ is $\frac{4h(\CO)\la_S(\CO)\mu_{\CO}}{\ka(\CO)}$ to one.
\end{proposition}
\prf
The question is, given $\CO\in O(S)$, how many conjugacy classes $[\ga]\in\CE_P^{p,\reg}(\Ga)$ are there such that $F_{\ga}\cap M(\Z)\cong\CO$.  Suppose $\CO$ is an order in the field $F$.  We saw in the proof of the previous proposition that given an embedding $\si:F\hra M(\Q)$ such that $\CO_{\si}=\CO$ and a fundamental unit $\ep\in\CO^{\x}$ the element $\ga=\si(\ep)$ satisfies both $F_{\ga}\cap M(\Z)\cong\CO$ and $[\ga]\in\CE_P^{p,\reg}(\Ga)$.  Conversely, it is clear that for $[\ga]\in\CE_P^{p,\reg}(\Ga)$, the element $\ga$ is a fundamental unit in the order $F_{\ga}\cap M(\Z)$.

Recall that $M(\Z)^{\x}$ acts on the set $\Si(\CO)$ from the left and that, by Proposition \ref{pro:embeddings}, the cardinality of the set $M(\Z)^{\x}\bs\Si(\CO)$ is $h(\CO)\la_S(\CO)$.  If $\ep,\d\in\CO^{\x}$ are fundamental units and $\si,\tau$ are embeddings of $F$ into $M(\Q)$ such that $\CO_{\si}=\CO_{\tau}=\CO$, then $\si(\ep)$ and $\tau(\d)$ are in the same $M(\Z)^{\x}$-conjugacy class if and only if there exists an automorphism $\al$ of $F$ such that $\al(\d)=\ep$ and the embeddings $\si$ and $\tau\circ\al$ are in the same class modulo $M(\Z)^{\x}$.

We will show that for each $\ga\in\CE_P^{p,\reg}(\Ga)$ the $M(\Z)^{\x}$-conjugacy class of $\ga$ decomposes into two $\Ga$-conjugacy classes.  If $\ga,\d\in\Ga$ are in the same $M(\Z)^{\x}$-conjugacy class we shall write $[[\ga]]=[[\d]]$.  For $\ga\in\Ga$, if $\ga$ is central in $M(\Q)$ then clearly the conjugacy class $[[\ga]]=\{\ga\}$.  Conversely, if $[[\ga]]=\{\ga\}$ then $\ga$ commutes with every element of $M(\Z)$ and, since $M(\Z)$ is an order in $M(\Q)$, it follows that $\ga$ commutes with every element of $M(\Q)$.  Hence the conjugacy class $[[\ga]]$ consists of a single element if and only if $\ga$ is central in $M(\Q)$.  Note that, by \cite{Reiner75}, Theorem 34.9, the image of $M(\Z)$ under the map $\det\!$ is equal to $\Z$, so there exist elements of $M(\Z)$ whose image under $\det\!$ is $-1$.

Let $\ga\in\CE_P^{p,\reg}(\Ga)$.  Then $\ga$ is not central in $M(\Q)$ so there exists $\d\neq\ga$ in $[[\ga]]$ such that $\ga=\eta^{-1}\d\eta$ for some $\eta\in M(\Z)^{\x}$ with $\det(\eta)=-1$.  Suppose further that $[\ga]=[\d]$, so that $\ga=\th^{-1}\d\th$ for some $\th\in\Ga$.  Then $\ga=\eta^{-1}\th\ga\th^{-1}\eta$, so $\xi=\eta^{-1}\th$ is an element of $F_\ga$ with $\det(\xi)=-1$.  By Lemma \ref{lem:unit} we then have $\xi\in F_{\ga}\cap M(\Z)^{\x}$, where $\CO=F_{\ga}\cap M(\Z)$ is an order in the totally complex quartic field $F_{\ga}$ and $\CO^{\x}=F_{\ga}\cap M(\Z)^{\x}$ is the group of units.  But since $F_{\ga}$ is totally complex quartic, Lemma \ref{lem:unit} tells us that $\det(\xi)=1$ and hence, by contradiction, $[\ga]\neq[\d]$.  This shows that each $M(\Z)^{\x}$-conjugacy class of elements in $\CE_P^{p,\reg}(\Ga)$ decomposes into at least two $\Ga$-conjugacy classes.

Suppose that $[[\ga]]=[[\d]]$ but $[\ga]\neq[\d]$.  Then there exists $\eta\in M(\Z)^{\x}$ such that $\ga=\eta^{-1}\d\eta$ and $\det(\eta)=-1$.  Let $\be\in\Ga$ also be such that $[[\ga]]=[[\be]]$ but $[\ga]\neq[\be]$.  Then there exists $\th\in M(\Z)^{\x}$ such that $\ga=\th^{-1}\be\th$ and $\det(\th)=-1$.  Then $\d=\eta\th^{-1}\be\th\eta^{-1}$, where $\det(\th\eta^{-1})=1$ so $[\d]=[\be]$.  Hence each $M(\Z)^{\x}$-conjugacy class of elements in $\CE_P^{p,\reg}(\Ga)$ decomposes into precisely two $\Ga$-conjugacy classes.

We conclude that the choice of a fundamental unit $\ep\in\CO^{\x}$ gives us $2h(\CO)\la_S(\CO)$ $\Ga$-conjugacy classes $[\ga]\in\CE_P^{p,\reg}(\Ga)$ such that $F_{\ga}\cap M(\Z)\cong\CO$.  We have seen that the choice of fundamental unit is determined up to inversion and multiplication by a root of unity so there are $2\mu_{\CO}$ choices for $\ep$.  We saw in the previous proposition that for a given embedding $\si$ of $F$ into $M(\Q)$ and choice of a fundamental unit $\ep\in\CO$ there is a conjugacy class $[\ga_{\si,\ep}]\in\CE_P^{p,\reg}(\Ga)$ with $\si(\ep)=\ga_{\si,\ep}$.  However, we saw above that if there exists an automorphism $\al$ of $F$ such that $\d=\al(\ep)$, then setting $\tau=\si\circ\al^{-1}$ we have $[\ga_{\si,\ep}]=[\ga_{\tau,\d}]$.  Such an automorphism exists if and only if $\d$ is a root of the minimal polynomial of $\ep$ over $\Q$.  The claim of the proposition folows.
\qed

\begin{lemma}
\label{lem:KaCondition}
Let $\ga\in\CE_P^{p,\rw}(\Ga)$ so that $\ga$ is conjugate in $G$ to $a_{\ga}b_{\ga}\in A^-B$ where
$$
a_{\ga}=\matrixfour{a}{a}{a^{-1}}{a^{-1}}\ \ \textrm{ and }\ \ b_{\ga}=\matrixtwo{R(\th)}{R(\phi)}
$$
for some $a\in(0,1)$, $\th,\phi\in\R$.  If $\ka(\CO_{\ga})>1$ then either $\th+\phi$ or $\th-\phi$ is in $\frac{\pi}{2}\Z\cup\frac{\pi}{3}\Z$.
\end{lemma}
\prf
The element $\ga$ is a fundamental unit in $\CO_{\ga}$.  We shall write $f(x)$ for its minimal polynomial.  The roots of $f(x)$ are the eigenvalues of $\ga$, which are $\al_1=ae^{i\th}$, $\al_2=ae^{-i\th}$, $\al_3=a^{-1}e^{i\phi}$, $\al_4=a^{-1}e^{-i\phi}$.

Suppose $\ka(\CO_{\ga})>1$.  Then there exist $i,j\in\{1,2,3,4\}$, $i\neq j$ and a root of unity $\ze$ such that $\al_i=\ze\al_j^{\pm 1}$.  If $\al_i=\ze\al_j$ then $\al_i\al_j^{-1}=\ze$ and either $\{i,j\}=\{1,2\}$ or $\{i,j\}=\{3,4\}$.  Without loss of generality take $i=1$ and $j=2$, then $\ze=\al_1\al_2^{-1}=e^{i2\th}$.  It follows that if $n\in\N$ is such that $\ze^n=1$ then $\al_1^{2n}=a^n\in(0,1)$ and so $F_{\ga}$ has a real quadratic subfield.  However, by Lemma \ref{lem:WeaklyNeat} this contradicts the assumption that $\ga\in\CE_P^{p,\rw}(\Ga)$, so we must have $\al_i=\ze\al_j^{-1}$.

By considering the possible values of $i$ and $j$ which would satisfy $\al_i=\ze\al_j^{-1}$ we see that either $\th+\phi$ or $\th+\phi$ is equal to $q\pi$, where $q\in\Q$ is such that $\ze=e^{iq\pi}$.  It remains to show what the possible values for $q$ are, that is, which roots of unity may occur in $\CO_{\ga}$.  We saw in the proof of Lemma \ref{lem:rootsOfUnity} that the only roots of unity which may occur in a quartic field are $\pm 1$ and roots of order $3,4,5,6,8,10$ or $12$.  If $\ze$ is a root of order $5,8,10$ or $12$ then $\ze$ generates a totally complex quartic field, however, $\Q(\ze)$ has the real quadratic subfield $\Q(\ze+\ze^{-1})$, which possibility is excluded in our case.  So we are left with the possibility that $\ze=\pm 1$ or $\ze$ is a root of order $3,4$ or $6$, which implies the claim of the lemma.
\qed

\begin{lemma}
Let $[\ga]\in\CE_P^{p,\reg}(\Ga)$.  Then
$$
\chi_1(\Ga_{\ga})=\rez{\mu_{\ga}}.
$$
\end{lemma}
\prf
By Proposition 1.10 of \cite{primgeoSL4},
$$
\chi_1(\Ga_{\ga})=\frac{\left[\Ga_{\ga,A}:\Ga_{\ga,A}'\right]}{\left[\Ga_{\ga}:\Ga_{\ga}'\right]},
$$
where $\Ga'\subset\Ga$ is a torsion free subgroup of finite index.  In the case under consideration we have that $\Ga_{\ga}$ is isomorphic as a group to the group of norm one elements in the order $\CO_{\ga}=F_{\ga}\cap M(\Z)$.  By Lemma \ref{lem:unit} this is equal to the group of units $\CO_{\ga}^{\x}$.  By Dirichlet's unit theorem and Lemma \ref{lem:rootsOfUnity} it follows that
$$
\Ga_{\ga}\cong\ep^{\Z}\x\mu(\ga),
$$
for some generator $\ep$ and where $\mu(\ga)$ is the finite cyclic group of roots of unity in $\CO_{\ga}^{\x}$.  It then follows that
$$
\Ga_{\ga}'\cong\ep^{k\Z},
$$
for some $k\in\N$.  Since every torsion element in $\Ga_{\ga}$ is sent to the identity under the projection onto $\Ga_{\ga,A}$ we also have the isomorphisms
$$
\Ga_{\ga,A}\cong\ep^{\Z}
$$
and
$$
\Ga_{\ga,A}'\cong\ep^{k\Z}.
$$
The lemma follows.
\qed

Let $[\ga]\in\CE_P^{p,\reg}(\Ga)$.  Then $\ga$ is conjugate in $G$ to a matrix $a_{\ga}b_{\ga}$, where $a_{\ga}=\diag(a,a,a^{-1},a^{-1})$ for some $0<a<1$ and
$$
b_{\ga}=\matrixtwo{R(\th)}{R(\phi)}
$$
for some $\th$, $\phi$.  Recall that the length $l_{\ga}$ of $\ga$ is defined to be $8|\log a|$ and $N(\ga)=e^{l_{\ga}}$.  For $\CO\in O(S)$ let $\ep_{\CO}$ be a fundamental unit of $\CO$ and $R(\CO)=2|\log|\ep_{\CO}||$ the regulator of $\CO$ and let $r(\CO)=e^{4R(\CO)}$.  Under the map $\th$, the element $\ga$ corresponds to a fundamental unit of the order $\th([\ga])$ and it is clear that
$$
r\left(\th\left([\ga]\right)\right)=N(\ga).
$$
We recall from the proof of Lemma \ref{lem:centralisers} that $\th([\ga])$ is an order in the field $\Q(ae^{i\th},a^{-1}e^{i\phi})$.  
By Lemma 2.10 of \cite{primgeoSL4} we have that $\tr\tilde{\si}(\ga)=4(1-\cos 2\th)(1-\cos 2\phi)$.  Let
$$
\nu(\ga)=\prod_{\al}\left(1-\frac{\al(\ga)}{|\al(\ga)|}\right),
$$
where the product is over the embeddings of $F_{\ga}$ into $\C$.  A simple calculation then shows that
$$
\tr\tilde{\si}(\ga)=(1-e^{i\th})(1-e^{-i\th})(1-e^{i\phi})(1-e^{-i\phi})=\nu(\ga).
$$

We summarise the results of this section in the following:
\begin{proposition}
\label{pro:regOrders}
The map $\th$ is surjective and $\frac{4h(\CO)\la_S(\CO)\mu_{\CO}}{\ka(\CO)}$ to one.  For $[\ga]\in\CE_P^{p,\reg}(\Ga)$ we have that $\chi_1(\Ga_{\ga})=1/\mu_{\ga}$, that $\tr\tilde{\si}(\ga)=\nu(\ga)$ and that $N(\ga)=r\left(\th\left([\ga]\right)\right)$.
\end{proposition}

\subsection{Class numbers of orders in totally complex quartic fields}
We are now in a position to prove our main theorem.
\begin{theorem}\textnormal{(Main Theorem)}
\label{thm:Main}
\index{Main Theorem}
Let $S$ be a finite, non-empty set of prime numbers with an even number of elements.  For $x>0$ let
$$
\pi_S(x)=\sum_{{\CO\in O^c(S)}\atop{R(\CO)\leq x}}\la_S(\CO)h(\CO),\index{$\pi_S(x)$}
$$
Then, as $x\ra\infty$ we have
$$
\pi_S(x)\sim\frac{e^{4x}}{8x}.
$$
\end{theorem}
\prf
By \cite{Borel69}, Proposition 17.6 we know that $\Ga$ has a weakly neat subgroup $\Ga'$ which is of finite index.  Hence there exists $n_{\Ga}\in\N$ such that all roots of unity which are eigenvalues of $\Ad(\ga)$ for some $\ga\in\Ga$ have order less than or equal to $n_{\Ga}$.  Define the sets
$$
B^0_1=\left\{\matrixtwo{R(\th)}{R(\phi)}\in B:\th,\phi\notin\pi\Z\cup\frac{\pi}{2}\Z\cup\cdots\cup\frac{\pi}{n_{\Ga}}\Z\right\}
$$
and
$$
B^0_2=B^0_1\cup\left\{\matrixtwo{R(\th)}{R(\phi)}\in B:(\th+\phi),(\th-\phi)\notin\frac{\pi}{2}\Z\cup\frac{\pi}{3}\Z\right\}.
$$
For $x>0$ define the functions
$$
\pi^{\rw}_1(x)=\sum_{{[\ga]\in\CE^{p,\rw}_P(\Ga)}\atop{N(\ga)\leq x}}\chi_1(\Ga_{\ga})\ \textrm{ and }\ \ 
\pi^{\rw}_2(x)=\sum_{{[\ga]\in\CE^{p,\rw}_P(\Ga);\ka(\CO_{\ga})=1}\atop{N(\ga)\leq x}}\chi_1(\Ga_{\ga}).
$$
Setting $B^0=B^0_1$ in Theorem 3.2 of \cite{primgeoSL4} we get as $x\ra\infty$
\begin{equation}
\label{eqn:Pi1Asymptotic}
\pi^{\rw}_1(x)\sim\frac{2x}{\log x}.
\end{equation}
Setting $B^0=B^0_2$ in Theorem 3.2 of \cite{primgeoSL4} we get, using Lemma \ref{lem:KaCondition}, as $x\ra\infty$
\begin{equation}
\label{eqn:Pi2Asymptotic}
\pi^{\rw}_2(x)\sim\frac{2x}{\log x}.
\end{equation}
We also define the following functions for $x>0$:
$$
\pi_{S,1}(x)=\sum_{{\CO\in O^c(S)}\atop{R(\CO)\leq x}}\frac{\la_S(\CO)}{\ka(\CO)}h(\CO)\ \textrm{ and }\ \ 
\pi_{S,2}(x)=\sum_{{\CO\in O^c(S);\ka(\CO_{\ga})=1}\atop{R(\CO)\leq x}}\la_S(\CO)h(\CO).
$$
By Proposition \ref{pro:regOrders} and the fact that $N(\ga)=e^{4R(\CO_{\ga})}$ we deduce from (\ref{eqn:Pi1Asymptotic}) that
\begin{equation}
\label{eqn:PiS1Asymptotic}
\pi_{S,1}(x)\sim\frac{e^{4x}}{8x}
\end{equation}
and from (\ref{eqn:Pi2Asymptotic}) that
\begin{equation}
\label{eqn:PiS2Asymptotic}
\pi_{S,2}(x)\sim\frac{e^{4x}}{8x}.
\end{equation}
For $x>0$ define
$$
\pi'_{S,2}(x)=\pi_{S,1}(x)-\pi_{S,2}(x)=\sum_{{\CO\in O^c(S);\ka(\CO_{\ga})>1}\atop{R(\CO)\leq x}}\frac{\la_S(\CO)}{\ka(\CO)}h(\CO).
$$
Then from (\ref{eqn:PiS1Asymptotic}) and (\ref{eqn:PiS2Asymptotic}) it follows that
$$
\pi'_{S,2}(x)=o\left(\frac{e^{4x}}{x}\right).
$$
Since $\ka(\CO)\leq 4$ for all $\CO\in O^c(S)$ we have
$$
\pi_{S,2}(x)\leq\pi_S(x)\leq\pi_{S,2}(x)+4\pi'_{S,2}(x)
$$
and the theorem follows.
\qed

We can also prove the following:

\begin{theorem}
\label{thm:Main2}
\index{Main Theorem}
Let $S$ be a finite, non-empty set of prime numbers with an even number of elements.  For $x>0$ let
$$
\tilde{\pi}_S(x)=\sum_{{\CO\in O^c(S)}\atop{R(\CO)\leq x}}\nu(\CO)\la_S(\CO)h(\CO),\index{$\tilde{\pi}_S(x)$}
$$
where
$$
\nu(\CO)=\rez{2\mu_{\CO}}\sum_{\ep}\prod_{\al}\left(1-\frac{\al(\ep)}{|\al(\ep)|}\right),\index{$\nu(\CO)$}
$$
where the sum is over the $2\mu_{\CO}$ different fundamental units of $\CO$ and the product is over the embeddings of $\CO$ into $\C$.

Then, as $x\ra\infty$ we have
$$
\pi_S(x)\sim\frac{e^{4x}}{2x}.
$$
\end{theorem}
\prf
Let $\CE_P^{p,\rnw}(\Ga)=\CE_P^{p,\reg}(\Ga)\smallsetminus\CE_P^{p,\rw}(\Ga)$ be the set of regular, non weakly neat elements in $\CE_P^p(\Ga)$.  Setting $B^0=B^{\reg}$ in Theorem 3.2 of \cite{primgeoSL4} we get
$$
\sum_{{[\ga]\in\CE^{p,\reg}_P(\Ga)}\atop{N(\ga)\leq x}}\chi_1(\Ga_{\ga})\sim\frac{2x}{\log x}
$$
and with (\ref{eqn:Pi1Asymptotic}) above it follows that
$$
\sum_{{[\ga]\in\CE^{p,\rnw}_P(\Ga)}\atop{N(\ga)\leq x}}\chi_1(\Ga_{\ga})=o\left(\frac{x}{\log x}\right).
$$
It follows from Lemma 2.10 of \cite{primgeoSL4} that $0\leq\tr\tilde{\si}(b_{\ga})\leq 16$ for all $\ga\in\CE_P(\Ga)$, and hence we have
$$
\sum_{{[\ga]\in\CE^{p,\rnw}_P(\Ga)}\atop{N(\ga)\leq x}}\chi_1(\Ga_{\ga})\tr\tilde{\si}(b_{\ga})=o\left(\frac{x}{\log x}\right).
$$
It then follows from Theorem 3.1 of \cite{primgeoSL4} that
$$
\sum_{{[\ga]\in\CE^{p,\rw}_P(\Ga)}\atop{N(\ga)\leq x}}\chi_1(\Ga_{\ga})\tr\tilde{\si}(b_{\ga})\sim\frac{8x}{\log x}.
$$
We can also deduce from (\ref{eqn:Pi1Asymptotic}) and (\ref{eqn:Pi2Asymptotic}) that
$$
\sum_{{[\ga]\in\CE^{p,\rw}_P(\Ga);\ka(\CO_{\ga})>1}\atop{N(\ga)\leq x}}\chi_1(\Ga_{\ga})\tr\tilde{\si}(b_{\ga})=o\left(\frac{x}{\log x}\right).
$$
The theorem then follows from Proposition \ref{pro:regOrders} using the same arguments as in the proof of the previous theorem.
\qed

The methods we have used to deduce the asymptotic result for class numbers from the Prime Geodesic Theorem were not sharp enough to preserve the error term which was proven there.  However, we make the following conjecture:

\begin{conjecture}
Under the conditions of Theorem \ref{thm:Main}, as $x\ra\infty$ we have
$$
\pi_S(x)=\rez{2}L(4x)+O\left(\frac{e^{3x}}{x}\right),
$$
where
$$
L(x)=\int_1^x\frac{e^t}{t}\,dt.
$$
\end{conjecture}

\end{document}